\documentclass[a4paper,10pt]{article}

\usepackage{amsfonts}
\usepackage{latexsym}
\usepackage{epsfig}
\usepackage{amssymb}
\usepackage{graphicx}
\usepackage{oldgerm}
\usepackage{theorem}

\def\ZV{{\bf Z}^{N}_{<}}
\def\ZW{{\bf N}^{N}_{<}}

\def\RV{{\bf R}^{N}_{<}}
\def\RW{{\bf R_+}^{N}_{<}}
\def\S{{\bf S}}
\def\u{{\bf u}}
\def\v{{\bf v}}
\def\x{{\bf x}}
\def\y{{\bf y}}
\def\z{{\bf z}}

\def\f{\widehat{f}_N}
\def\h{\widehat{h}_N}

\def\gT{\widehat{g}_N^T}
\def\p{\widehat{p}_N}
\def\cN{{\cal N}_N}
\def\hN{\widehat{\cal N}_N}

\theorembodyfont{\itshape}

\newtheorem{thm}{Theorem}[section]
\newtheorem{lem}[thm]{Lemma}
\newtheorem{cor}[thm]{Corollary}

\newcommand{\SSC}[1]{\section{#1}\setcounter{equation}{0}}
\newcommand{\qed}{\hbox{\rule[-2pt]{3pt}{6pt}}}

\begin{document}
\begin{center}
{\bf \Large{Functional central limit theorems for vicious walkers}}
\vskip 3mm
{\sc Makoto Katori and Hideki Tanemura}
\\
{\it Chuo University} and {\it Chiba University} 
\vskip 3mm
Dedicated to Professor Tokuzo Shiga on his 60th birthday.
\end{center}

\pagestyle{plain}
\vskip 0.3cm

\begin{small}
\noindent  
{\bf Abstract.} 
We consider the diffusion scaling limit of the vicious walker model
that is a system of nonintersecting random walks.
We prove a functional central limit theorem for the model
and derive two types of nonintersecting Brownian motions, 
in which the nonintersecting condition is imposed
in a finite time interval $(0,T]$ for the first type and
in an infinite time interval $(0,\infty)$
for the second type, respectively.
The limit process of the first type is a temporally inhomogeneous 
diffusion, and that of the second type is a temporally homogeneous 
diffusion that is identified with a Dyson's model of Brownian 
motions studied in the random matrix theory. 
We show that these two types of processes are related to
each other by a multi-dimensional generalization of
Imhof's relation, whose original form relates the
Brownian meander and the three-dimensional Bessel process.
We also study the vicious walkers with wall restriction
and prove a functional central limit theorem in the diffusion 
scaling limit.

\vskip 0.5cm



\end{small}

\normalsize

\SSC{Introduction}

The system of one-dimensional symmetric simple random walks, 
in which none of walkers have met others in a given time period, 
is called the vicious walker model. 
( See Fisher's paper \cite{Fis84}, in which it was
introduced as a model of statistical mechanics.)
The purpose of this paper is to study the scaling limit
of vicious walkers as a stochastic process.
Since each random walk tends to a Brownian motion
in the diffusion scaling limit, an interacting system of
$N$ Brownian motions will be constructed as the scaling limit of
vicious walkers with an arbitrary finite number of walkers $N$.
We show that a functional central limit theorem for vicious walkers holds
and the limit process ${\bf X}(t)=(X_1(t),X_2(t),\dots,X_N(t))$
is a temporally inhomogeneous diffusion, that is,
its transition probability depends on the time interval $(0,T]$
in which the nonintersecting condition is imposed.
We claim that when $N=2$, 
the process  $(X_2(t)-X_1(t))/\sqrt{2}$ 
is a one-dimensional Brownian motion conditioned to stay positive 
during a finite time interval $(0,T]$, 
which is called a {\it Brownian meander} in \cite{RY98,Yor92}.

We also study the temporally homogeneous diffusion process 
${\bf Y}(t)=(Y_1(t),Y_2(t),\dots,Y_N(t))$,
which is obtained from the previous process by taking $T \to \infty$.
The process ${\bf Y}(t)$ is the {\it Doob h-transform} \cite{Do84}
of the absorbing Brownian motion in a Weyl chamber
with harmonic function $h_{N}(\x)= \prod_{1\le i<j\le N}(x_j-x_i)$,
and can be regarded as a system of Brownian motions 
with the drift terms acting as the repulsive two-body forces 
proportional to the inverse of distances between particles \cite{Gra99}.
In other words, if we set $T \to \infty$,
the scaling limit of vicious walkers can
realize a Dyson's Brownian motion model 
studied in the random matrix theory \cite{Dys62,Meh91}.
We show the following relation between 
the processes ${\bf X}(t)$ and ${\bf Y}(t)$:
$$
P({\bf X}(\cdot)\in dw)
=\overline{c}_N (T) P({\bf Y}(\cdot)\in dw)\frac{1}{h(w(T))},
$$
where $\overline{c}_N(T)$ is the normalization constant.
This equality is a generalization of Imhof's formula,
which relates the Brownian meander and the three-dimensional 
Bessel process \cite{Imh84}.

The Gaussian ensembles of random matrices can be regarded
as the thermodynamical equilibrium of Coulomb gas system
and that is the reason why Dyson introduced a one-dimensional
model of interacting Brownian particles with 
(two-dimensional) Coulomb repulsive potentials \cite{Dys62,Meh91}.
Similar relations between our processes and random matrix ensembles 
can be seen. The distribution of ${\bf Y}(t)$ is described by using
the probability density of eigenvalues
of random matrices in the Gaussian unitary ensemble (GUE)
with variance $t$.
Pandey and Mehta \cite{MP83, PM83} introduced 
a Gaussian ensemble of Hermitian matrices depending 
on a parameter $\alpha\in [0,1]$.
When $\alpha =0$, the ensemble is 
the Gaussian orthogonal ensemble (GOE),
and when $\alpha =1$, it is the GUE.
We will find that the probability density function of 
$\sqrt{T/\{t(2T-t)\}} {\bf X}(t)$ coincides with that
of eigenvalues of matrices in the Pandey-Mehta ensemble
with $\alpha = \sqrt{(T-t)/T}$.

We also study vicious walkers with wall restriction
and prove the functional central limit theorem in the diffusion 
scaling limit.
In this case the obtained temporally inhomogeneous diffusion process
is a system of nonintersecting Brownian meanders, 
which is related to the nonstandard classes of 
random matrices \cite{KTNK*1,KO01}.

\SSC{Statement of Results}
\subsection{Vicious walkers without wall restriction}

Let $(\{\S_j\}_{j\ge 0},P^{\z})$ be the $N$-dimensional Markov chain
starting from $\z=(z_1,z_2,\dots, z_N)$, 
such that the coordinates $S_j^k, k=1,2,\dots,N$,
are independent simple random walks on ${\bf Z}$.
We always take the starting point $\z$ from the set
\begin{equation}
\ZV =
\{ \z =(z_1,z_2,\dots,z_N) \in (2{\bf Z})^N \; ; \; 
z_{k+1}-z_k \in 2{\bf Z_+}, \; k=1,\ldots, N-1\}, 
\nonumber
\end{equation} 
where ${\bf Z_+}$ is the set of positive integers.
Now we consider the condition that any of walkers 
does not meet other walkers up to time $m$, i.e.
\begin{equation}
S_j^1 < S_j^2 < \cdots < S_j^N,
\quad 0\le j \le m.
\label{eqn:nonint}
\end{equation}
We dente by $Q^{\z}_m$ the conditional probability of $P^{\z}$
under the event 
$\Lambda_m = \{S_j^1 < S_j^2 < \cdots < S_j^N, \; 0\le j \le m \}$.
The process $(\{\S_j\}_{j\ge 0},Q^{\z}_m)$ is called 
the vicious walkers (up to time $m$) (see Fisher \cite{Fis84}).
\begin{figure}[ht]
\begin{center}
\unitlength 0.1in
\begin{picture}( 36.0000, 22.6500)(  6.0000,-28.6500)
%
\special{pn 13}%
\special{pa 2000 2800}%
\special{pa 2200 2600}%
\special{pa 2000 2400}%
\special{pa 2200 2200}%
\special{pa 2000 2000}%
\special{pa 1800 1800}%
\special{pa 2000 1600}%
\special{pa 1800 1400}%
\special{pa 1600 1200}%
\special{pa 1800 1000}%
\special{pa 1800 1000}%
\special{pa 1600 800}%
\special{pa 1600 800}%
\special{pa 1600 800}%
\special{fp}%
%
\special{pn 13}%
\special{pa 1600 2800}%
\special{pa 1400 2600}%
\special{pa 1200 2400}%
\special{pa 1400 2200}%
\special{pa 1600 2000}%
\special{pa 1400 1800}%
\special{pa 1200 1600}%
\special{pa 1400 1400}%
\special{pa 1200 1200}%
\special{pa 1000 1000}%
\special{pa 1000 1000}%
\special{pa 1000 1000}%
\special{pa 800 800}%
\special{pa 800 800}%
\special{fp}%
%
\special{pn 13}%
\special{pa 2400 2800}%
\special{pa 2600 2600}%
\special{pa 2800 2400}%
\special{pa 2600 2200}%
\special{pa 2800 2000}%
\special{pa 2600 1800}%
\special{pa 2400 1600}%
\special{pa 2200 1400}%
\special{pa 2400 1200}%
\special{pa 2600 1000}%
\special{pa 2600 1000}%
\special{pa 2400 800}%
\special{pa 2400 800}%
\special{pa 2400 800}%
\special{fp}%
%
\special{pn 13}%
\special{pa 3200 2800}%
\special{pa 3000 2600}%
\special{pa 3200 2400}%
\special{pa 3400 2200}%
\special{pa 3200 2000}%
\special{pa 3400 1800}%
\special{pa 3200 1600}%
\special{pa 3400 1400}%
\special{pa 3200 1200}%
\special{pa 3000 1000}%
\special{pa 3000 1000}%
\special{pa 3200 800}%
\special{pa 3200 800}%
\special{fp}%
%
\special{pn 8}%
\special{pa 600 2800}%
\special{pa 4200 2800}%
\special{fp}%
%
\special{pn 8}%
\special{pa 2000 2800}%
\special{pa 2000 600}%
\special{fp}%
\special{sh 1}%
\special{pa 2000 600}%
\special{pa 1980 668}%
\special{pa 2000 654}%
\special{pa 2020 668}%
\special{pa 2000 600}%
\special{fp}%
%
\special{pn 8}%
\special{pa 600 800}%
\special{pa 3600 800}%
\special{dt 0.045}%
\put(16.0000,-29.5000){\makebox(0,0){$z_1$}}%
\put(20.0000,-29.5000){\makebox(0,0){$z_2$=0}}%
\put(24.0000,-29.5000){\makebox(0,0){$z_3$}}%
\put(32.0000,-29.5000){\makebox(0,0){$z_4$}}%
\put(15.3000,-18.0000){\makebox(0,0){$S^1$}}%
\put(19.5000,-18.0000){\makebox(0,0){$S^2$}}%
\put(27.6000,-18.0000){\makebox(0,0){$S^3$}}%
\put(35.1000,-18.0000){\makebox(0,0){$S^4$}}%
\put(22.4000,-7.2000){\makebox(0,0){\it m=10}}%
\end{picture}%

\end{center}
\caption{{\it An example of vicious walks without wall restriction.}}
\label{fig1}
\end{figure}
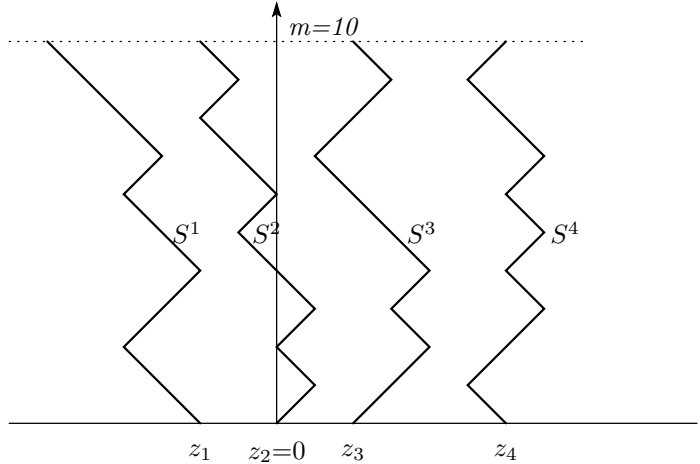

For $T>0$ and $\z\in \ZV$, 
we consider probability measures $\mu_{L,T}^{\z}$, $L\ge 1$,
on the space of continuous paths $C( [0,T]\to{\bf R}^N)$ 
defined by
\begin{equation}
\mu_{L,T}^{\z}(\cdot) 
= Q_{L^2T}^{\z} \left(\frac{1}{{L}}\S(L^2t) \in \cdot \right),
\end{equation}
where $\S(t), t\ge 0$, is the interpolation of the random walk
$\S_j, j=0,1,2,\dots$.
We study the limit distribution of the probability
$\mu_{L,T}^{\z}$, $L\to \infty$.

We put $\RV=\{\x\in {\bf R}^N : x_1 < x_2 < \cdots < x_N \}$,
which is called the Weyl chamber \cite{FH91}.
By virtue of the Karlin-McGregor formula \cite{KM59_1,KM59_2}, 
the transition density function  $f_N (t, \y|\x)$ of
the absorbing Brownian motion in $\RV$
and the probability $\cN (t, \x)$ that the Brownian motion 
started at $\x\in\RV$ does not
hit the boundary of $\RV$ up to time $t>0$
are given by
\begin{equation}
f_{N}(t, \y|\x)= 
\det_{1 \leq i, j \leq N}
\left( (2 \pi t)^{-1/2} \ e^{-(x_{j}-y_{i})^2/2t} \right),
\: \x,\y \in \RV,
\end{equation}
and
\begin{equation}
\cN (t, \x ) = 
\int_{\RV} d\y f_N (t, \y |\x ).
\nonumber
\end{equation}

For an even integer $n$ and an antisymmetric $n\times n$ matrix 
$A = (a_{ij})$ we put
\begin{equation}
{\rm Pf}_{1\le i<j \le n}(a_{ij})
= \frac{1}{(n/2)!} \sum_{\sigma} sgn (\sigma)
a_{\sigma(1)\sigma(2)}a_{\sigma(3)\sigma(4)}
\cdots a_{\sigma(n-1)\sigma(n)},
\end{equation}
where the summation is extended over all permutations $\sigma$
of $(1,2,\dots,n)$ with restriction
$\sigma(2k-1)<\sigma(2k)$, $k=1,2,\dots,n/2$.
This expression is known as the {\it Pfaffian} 
(see Stembridge \cite{Stem90}).
Then we have the following lemma,
which is a consequence of the identity given by
de Bruijn \cite{deBr55} as shown in Section 4.

\vskip 3mm
\begin{lem}
\label{lem:Nor_V}
For $t > 0$, $\x \in \RV$,
\begin{eqnarray}
\cN (t,\x) &=&
\left\{
   \begin{array}{ll}
      {\rm Pf}_{1\le i<j\le N}
      \displaystyle{F_{ij}(t,\x)},
      & \mbox{if} \ N=\mbox{even} 
\\
      {\rm Pf}_{1\le i<j\le N+1}
      \displaystyle{F_{ij}(t,\x)},
      & \mbox{if} \ N=\mbox{odd}, \\
   \end{array}
\right.
\end{eqnarray}
where
\begin{eqnarray}
F_{ij}(t,\x) &=&
\left\{
   \begin{array}{ll}
\Psi \left( \displaystyle{\frac{x_j-x_i}{2\sqrt{t}}} \right),
      & \mbox{if} \ 1\le i, j \le N,
\\
1, & \mbox{if} \ 1\le i \le N, j=N+1,
\\
-1, & \mbox{if} \ i=N+1, 1\le j \le N,
\\
0, & \mbox{if} \ i=N+1, j=N+1,
\\
   \end{array}
\right.
\end{eqnarray}
and $\Psi(u)= (2/ \sqrt{\pi}) \int_0^u e^{-v^2}dv$.
\end{lem}

We put 
\begin{equation}
h_N(\x) = \prod_{1\le i <j \le N}( x_j - x_i).
\end{equation}
The first main result is the following theorem.
\begin{thm}
\label{thm:VT}
{\rm (i)}
For any fixed $\z \in \ZV$ and $T>0$, as $L\to \infty$,
$\mu_{L,T}^{\z}(\cdot)$ converges weakly to the law of 
the temporally inhomogeneous diffusion process
${\bf X}(t)=(X_1(t),X_2(t),\dots,X_N(t)), t\in [0,T]$,
with transition density $g^T_N(s,\x,t,\y)$;
\begin{equation}
\label{eqn:gn0}
g_{N}^{T}(0, {\bf 0}, t, \y)
=c_{N}T^{N(N-1)/4}t^{-N^2/2}  
\exp\left\{ -\frac{|\y|^2}{2t} \right\}
h_N(\y)\cN (T-t,\y),
\end{equation}
\begin{equation}
g_N^T(s,\x, t, \y )
= \frac{f_{N}(t-s, \y|\x)\cN (T-t,\y)}{\cN (T-s,\x)},
\end{equation}
for $0 \leq s < t \le T,\; \x, \y \in \RV,$
where 
$c_N = 2^{-N/2}/
\prod_{j=1}^N \Gamma(j/2)$.

\noindent
{\rm (ii)}
The diffusion process ${\bf X}(t)$
solves the following equation:
\begin{equation}
\label{eq:SDE1}
X_i (t) = B_i(t) + \int_0^t b_i^T(s, {\bf X}(s)) ds,
\quad t\in [0,T],\quad i=1,2,\dots, N,
\end{equation}
where $B_i(t)$, $i=1,2,\dots,N,$ 
are independent one-dimensional Brownian motions and
$$
b_i^T (t,\x) = \frac{\partial}{\partial x_i}\ln \cN (T-t, \x),
\quad i=1,2,\dots,N.
$$
\end{thm}

Although $\mu_{L,T}^{\z}(\cdot)$ is the probability measure
defined on $C([0,T]\to {\bf R}^N)$,
it can be regarded as that on $C([0,\infty)\to {\bf R}^N)$
concentrated on the set
$$
\{ w \in C([0,\infty)\to {\bf R}^N) : w(t)= w(T),\: t \ge T\}.
$$
Next we consider the case that $T=T_L$ goes to infinity 
as $L \to \infty$.

\begin{cor}
\label{cor:VI}
{\rm (i)}
Let $T_L$ be an increasing function of $L$ with
$T_L \to \infty$ as $L\to\infty$.
For any fixed $\z \in \ZV$, as $L\to \infty$,
$\mu_{L,T_L}^{\z}(\cdot)$ converges weakly to the law of 
the temporally homogeneous diffusion process
${\bf Y}(t)=(Y_1(t),Y_2(t),\dots,Y_N(t)), t\in [0,\infty)$,
with transition density $p_N(s,\x,t,\y)$;
\begin{equation}
p_{N}(0, {\bf 0}, t, \y)
= c^{\prime}_{N}t^{-N^2/2} 
\exp\left\{ -\frac{|\y|^2}{2t} \right\} h_N(\y)^2,
\label{eqn:pn0}
\end{equation}
\begin{equation}
p_N(s,\x, t, \y )
= \frac{1}{h_N(\x)}f_{N}(t-s, \y|\x)h_N(\y),
\end{equation}
for $0 \leq s < t < \infty,\; \x, \y \in \RV,$
where 
$c^{\prime}_{N}=(2\pi)^{-N/2}
/\prod_{j=1}^{N} \Gamma(j)$.

\noindent
{\rm (ii)}
The diffusion process ${\bf Y}(t)$ solves 
the equation of Dyson's Brownian motion model \cite{Dys62} :
\begin{equation}
\label{eq:SDE2}
Y_i (t) = B_i(t) 
+ \sum_{1\le j \le N,j\not=i}\int_0^t \frac{1}{Y_i(s)-Y_j(s)}ds,
\quad t\in [0,\infty), \quad i=1,2,\dots, N.
\end{equation}
\end{cor}

\subsection{Vicious walkers with wall restriction}

In this subsection, we impose the condition
\begin{equation}
S_j^1  \geq 0, \quad
1 \le j \le m,
\label{eqn:wall}
\end{equation}
in addition to (\ref{eqn:nonint})
and take the starting point $\z$ from the set
\begin{equation}
\ZW =
\{ \z =(z_1,z_2,\dots,z_N) \in (2{\bf N})^N \; ; \; 
z_{k+1}-z_k \in 2{\bf Z_+}, \; k=1,\ldots, N-1\}, 
\nonumber
\end{equation} 
where ${\bf N}$ is the set of non-negative integers.
That is, there assumed to be a wall at the origin
and all walkers can walk only in the region $[0,\infty)$. 
We dente by $\widehat{Q}^{\z}_m$ the conditional probability of 
$P^{\z}$ under the event 
$\widehat{\Lambda}_m = \{0\le S_j^1 < S_j^2 < \cdots < S_j^N,
\; 0\le j \le m \}$.
The process $(\{{\bf S}_j\}_{j\ge 0},\widehat{Q}^{\z}_m)$ is regarded
as the vicious walkers with wall restriction (up to time $m$)
\cite{KGV00}. 
\begin{figure}[h]
\begin{center}
\unitlength 0.1in
\begin{picture}( 40.0000, 22.6500)(  6.0000,-28.6500)
%
\special{pn 8}%
\special{pa 600 2800}%
\special{pa 4200 2800}%
\special{fp}%
%
\special{pn 8}%
\special{pa 2000 2800}%
\special{pa 2000 600}%
\special{fp}%
\special{sh 1}%
\special{pa 2000 600}%
\special{pa 1980 668}%
\special{pa 2000 654}%
\special{pa 2020 668}%
\special{pa 2000 600}%
\special{fp}%
%
\special{pn 8}%
\special{pa 1600 786}%
\special{pa 4600 786}%
\special{dt 0.045}%
\put(20.0000,-29.5000){\makebox(0,0){$z_1$=0}}%
\put(24.0000,-29.5000){\makebox(0,0){$z_2$}}%
\put(32.0000,-29.5000){\makebox(0,0){$z_3$}}%
\put(36.0000,-29.5000){\makebox(0,0){$z_4$}}%
\put(24.0000,-18.0000){\makebox(0,0){$S^1$}}%
\put(31.2000,-18.0000){\makebox(0,0){$S^2$}}%
\put(35.5000,-18.0000){\makebox(0,0){$S^3$}}%
\put(39.8000,-18.0000){\makebox(0,0){$S^4$}}%
\put(22.1000,-7.1000){\makebox(0,0){\it m=10}}%
%
\special{pn 13}%
\special{pa 2000 2800}%
\special{pa 2200 2600}%
\special{pa 2000 2400}%
\special{pa 2200 2200}%
\special{pa 2400 2000}%
\special{pa 2200 1800}%
\special{pa 2400 1600}%
\special{pa 2200 1400}%
\special{pa 2400 1200}%
\special{pa 2600 1000}%
\special{pa 2600 1000}%
\special{pa 2400 800}%
\special{pa 2400 800}%
\special{pa 2400 800}%
\special{fp}%
%
\special{pn 13}%
\special{pa 2400 2800}%
\special{pa 2600 2600}%
\special{pa 2800 2400}%
\special{pa 2600 2200}%
\special{pa 2800 2000}%
\special{pa 3000 1800}%
\special{pa 2800 1600}%
\special{pa 2600 1400}%
\special{pa 2800 1200}%
\special{pa 3000 1000}%
\special{pa 3000 1000}%
\special{pa 3200 800}%
\special{pa 3200 800}%
\special{pa 3200 800}%
\special{fp}%
%
\special{pn 13}%
\special{pa 3200 2800}%
\special{pa 3000 2600}%
\special{pa 3200 2400}%
\special{pa 3400 2200}%
\special{pa 3200 2000}%
\special{pa 3400 1800}%
\special{pa 3600 1600}%
\special{pa 3400 1400}%
\special{pa 3600 1200}%
\special{pa 3400 1000}%
\special{pa 3400 1000}%
\special{pa 3600 800}%
\special{pa 3600 800}%
\special{pa 3600 800}%
\special{fp}%
%
\special{pn 13}%
\special{pa 3600 2800}%
\special{pa 3400 2600}%
\special{pa 3600 2400}%
\special{pa 3800 2200}%
\special{pa 3600 2000}%
\special{pa 3800 1800}%
\special{pa 4000 1600}%
\special{pa 4200 1400}%
\special{pa 4000 1200}%
\special{pa 4200 1000}%
\special{pa 4000 800}%
\special{pa 4000 800}%
\special{fp}%
\end{picture}%

\end{center}
\caption{
{\it An example of vicious walks with wall restriction.}}
\label{fig2}
\end{figure}
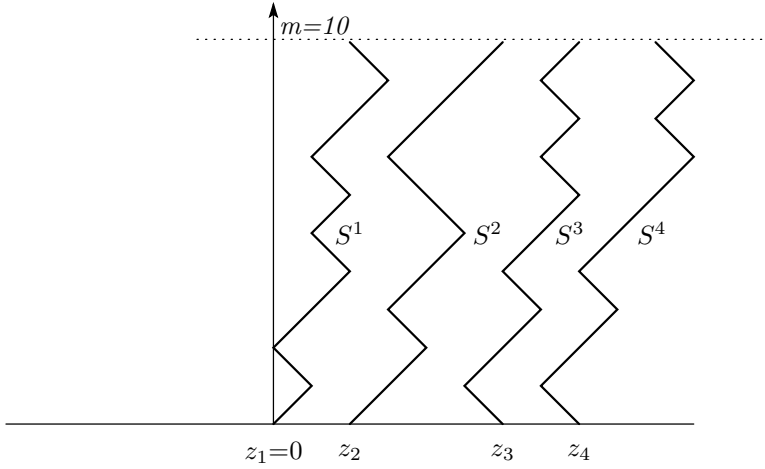

For $T>0$ and $\z\in \ZW$,  we consider probability measures 
$\widehat{\mu}_{L,T}^{\z}$, $L\ge 1$,
on $C( [0,T]\to{\bf R}^n)$ defined by
\begin{equation}
\widehat{\mu}_{L,T}^{\z}(\cdot) 
= \widehat{Q}_{L^2T}^{\z} \left(\frac{1}{{L}}\S(L^2t) \in \cdot \right).
\end{equation}
We study the limit distribution of the probability
$\widehat{\mu}_{L,T}^{\z}$, $L\to \infty$.

We put 
$\RW = \{\x\in {\bf R}^N : 0\le x_1 < x_2 < \cdots < x_N \}$. 
Then the transition density function  $\f (t, \y|\x)$ of
the absorbing Brownian motion in $\RW$
and the probability $\hN (t, \x)$ that the Brownian motion 
started at $\x\in \RW$ does not
hit the boundary of $\RW$ up to time $t>0$
are given by
\begin{equation}
\f (t, \y|\x)= 
\det_{1 \leq i, j \leq N}
\left( (2 \pi t)^{-1/2} \ 
( e^{-(x_{j}- y_{i})^2/2t}
- e^{-(x_{j}+ y_{i})^2/2t})
\right),
\: \x,\y\in \RW,
\end{equation}
and
\begin{equation}
\hN (t,\x)= \int_{\RW}d\y \f (t, \y|\x).
\end{equation}
Then we have the following lemma, which is proved in Section 4
as a consequence of the identity given by de Bruijn \cite{deBr55}.
\begin{lem}
\label{lem:Nor_W}
For $t>0$, $x\in \RW$,
\begin{eqnarray}
\hN (t,\x) &=&
\left\{
   \begin{array}{ll}
      {\rm Pf}_{1\le i<j\le N}
      \displaystyle{\widehat{F}_{ij}(t,\x)},
      & \mbox{if} \ N=\mbox{even},
\\
      {\rm Pf}_{1\le i<j\le N+1}
      \displaystyle{\widehat{F}_{ij}(t,\x)},
      & \mbox{if} \ N=\mbox{odd}, \\
   \end{array}
\right.
\label{hNtx}
\end{eqnarray}
where
\begin{eqnarray}
\widehat{F}_{ij}(t,\x) &=&
\left\{
   \begin{array}{ll}
\widehat{\Psi}\left(
  \displaystyle{\frac{x_i}{\sqrt{2t}}, \frac{x_j}{\sqrt{2t}}} \right),
      & \mbox{if} \ 1\le i, j \le N,
\\
\Psi \left( \displaystyle{\frac{x_i}{\sqrt{2t}}}\right),
 & \mbox{if} \ 1\le i \le N, j=N+1,
\\
-\Psi \left( \displaystyle{\frac{x_j}{\sqrt{2t}}}\right),
& \mbox{if} \ i=N+1, 1\le j \le N,
\\
0, & \mbox{if} \ i=N+1, j=N+1,
\\
   \end{array}
\right.
\label{hF}
\end{eqnarray}
and
\begin{eqnarray}
\widehat{\Psi}(u_1,u_2)
&=& \frac{2}{\pi}
\left[ \int_0^{u_1} dv_1 \int_{u_1 -u_2}^{u_2 -u_1}dv_2
\exp \{ -v_1^2 - (v_1-v_2)^2 \}
\right.
\nonumber
\\
&&\quad - \left. 
\int_{u_1}^{u_2} dv_1 \int_{u_2 -u_1}^{u_1 +u_2}dv_2
\exp \{ -v_1^2 - (v_1-v_2)^2 \}
\right].
\label{hIt}
\end{eqnarray}
\end{lem}

We put
\begin{equation}
\h (\x )=\prod_{1 \leq i < j \leq N}(x_{j}^{2}-x_{i}^{2}) 
\prod_{i=1}^{N} x_{i}.
\end{equation}
Then we can obtain the following result.
\vskip 3mm
\begin{thm}
\label{thm:WT}
{\rm (i)}
For any fixed $\z \in \ZW$ and $T>0$, as $L\to \infty$,
$\widehat{\mu}_{L,T}^{\z}(\cdot)$ converges weakly to the law of 
the temporally inhomogeneous diffusion process
$\widehat{\bf X}(t)=(\widehat{X}_1(t), \widehat{X}_2(t),
\dots, \widehat{X}_N(t))$, $t\in [0,T]$,
with transition density $\gT (s,\x,t,\y)$;
\begin{equation}
\gT (0, {\bf 0}, t, \y)
=\widehat{c}_{N}T^{N^2/2}t^{-N(2N+1)/2} 
\exp\left\{ -\frac{|\y|^2}{2t} \right\} \h(\y)
\hN (T-t,\y),
\end{equation}
\begin{equation}
\gT (s,\x, t, \y )
= \frac{\f (t-s, \y|\x)\hN (T-t,\y)}{\hN(T-s,\x)}, 
\end{equation}
\end{thm}
for $0 \leq s < t \le T, \x,\y \in \RW$,
where 
$\widehat{c}_N = 1/\prod_{j=1}^N \Gamma(j)$.

\noindent
{\rm (ii)}
The diffusion process $\widehat{\bf X}(t)$
solves the following equation:
\begin{equation}
\widehat{X}_i (t) = B_i(t) 
+ \int_0^t \widehat{b}_i^T(s, \widehat{X}(s)) ds,
\quad t\in [0,T],\quad i=1,2,\dots,N,
\end{equation}
where
$$
\widehat{b}_i^T (t,\x) = \frac{\partial}{\partial x_i}\ln \hN (T-t, \x),
\quad i=1,2,\dots,N.
$$
\vskip 3mm

Next we consider the case that $T=T_L$ goes to infinity 
as $L \to \infty$.

\vskip 3mm
\begin{cor}
\label{cor:WI}
{\rm (i)}
Let $T_L$ be an increase function of $L$ with $T_L\to\infty$
as $L\to\infty$.
For any fixed $\z \in \ZW$, as $L\to \infty$,
$\widehat{\mu}_{L,T_L}^{\z}(\cdot)$ converges weakly to the law of 
the temporally homogeneous diffusion process
$\widehat{\bf Y}(t)=(\widehat{Y}_1(t), \widehat{Y}_2(t),
\dots, \widehat{Y}_N(t))$, $t\in [0,\infty)$,
with transition density $\p (s,\x,t,\y)$;
\begin{equation}
\p (0, {\bf 0}, t, \y)
= \widehat{c^{\prime}}_N
t^{-N(2N+1)/2} \exp\left\{-\frac{|\y|^2}{2t}\right\}\h(\y)^2,
\label{eqn:pnh0}
\end{equation}
\begin{equation}
\p (s,\x, t, \y )
= \frac{1}{\h(\x)}\f (t-s, \y | \x )\h(\y),
\end{equation}
for $0\le s < t < \infty$, $ \x, \y \in \RW$,
where 
$\widehat{c^{\prime}}_{N}=(2/\pi)^{N/2}
/\prod_{j=1}^N \Gamma(2j)$.

\noindent
{\rm (ii)}
The diffusion process $\widehat{\bf Y}(t)$
solves the following equation:
\begin{equation}
\widehat{Y}_i (t) = B_i(t) + 
\int_0^t \frac{1}{\widehat{Y}_i(s)}ds
+ \sum_{1\le j \le N,j\not=i} \left\{
\int_0^t \frac{1}{\widehat{Y}_i(s)-\widehat{Y}_j(s)} ds
+\int_0^t \frac{1}{\widehat{Y}_i(s)+\widehat{Y}_j(s)} ds
\right\},
\end{equation}
$t\in [0,\infty), i=1,2,\dots,N$.
\end{cor}
\vskip 3mm
\subsection{Remarks}

\noindent (i)
The process ${\bf X}(t)$ ( resp. $\widehat{\bf X}(t)$ )
represents the system of $N$ Brownian motions
( resp. $N$ Brownian meanders ) started from the origin
conditioned not to collide up to time $T$.
A limit theorem for one-dimensional random walk 
conditioned to stay positive
was firstly observed by Spitzer \cite{Spiz87}
and then studied and generalized by many probabilists
\cite{B72,Idl74,Bol76,Dur78}.
For two-dimensional random walk conditioned to stay in a cone,
a limit theorem was proved by Shimura \cite{Shi91}.
Our theorems are multi-dimensional versions of these limit theorems.

\noindent (ii)
The process ${\bf Y}(t)$ ( resp. $\widehat{\bf Y}(t)$ )
represents the system of $N$ Brownian motions 
( resp. $N$ three-dimensional Bessel processes )
conditioned never to collide.
The function $h_{N}(\x)$ ( resp. $\h (\x)$ )
is a strictly positive harmonic function for
the absorbing Brownian motions in the Weyl chamber $\RV$ (resp. $\RW$ ).
The process ${\bf Y}(t)$ ( resp. $\widehat{\bf Y}(t)$ )
is the corresponding {\it Doob h-transform} \cite{Do84,Gra99}.
A functional central limit theorem to the process ${\bf Y}(t)$
was also discussed in a recent paper by 
O'Connell and Yor \cite{OY02}.

\noindent(iii)
The relation between the Brownian meander and the three-dimensional 
Bessel process was discussed in Imhof \cite{Imh84}.
From our results Imhof's relation is generalized as follows:

\noindent
({\bf Without wall restriction})
For any $t_0 =0 < t_1 < \cdots < t_\ell =T$, $\ell \in {\bf Z_+}$,
\begin{equation}
\prod_{i=1}^{\ell} g_N^T (t_{i-1},\y_{i-1}, t_{i}, \y_i)
= \overline{c}_N T^{N(N-1)/4}
\prod_{i=1}^{\ell} p_N (t_{i-1},\y_{i-1}, t_{i}, \y_i)
\frac{1}{h_N(\y_{\ell})},
\end{equation}
for any $\y_{i}\in \RV, i=1,2,\dots, \ell$,
where $\y_0 = {\bf 0}$ and
$$
\overline{c}_N = 
\frac{c_N}{c^{\prime}_{N}}
=\pi^{N/2}\prod_{j=1}^N \frac{\Gamma (j)}{\Gamma (j/2)}.
$$

\noindent
({\bf With wall restriction})
For any $t_0 =0 < t_1 < \cdots < t_\ell =T$, $\ell \in {\bf Z_+}$
\begin{equation}
\prod_{i=1}^{\ell} \gT (t_{i-1},\y_{i-1}, t_{i}, \y_i)
= \widetilde{c}_N T^{N^2/2} 
\prod_{i=1}^{\ell} \p (t_{i-1},\y_{i-1}, t_{i}, \y_i)
\frac{1}{\h(\y_{\ell})},
\end{equation}
for any $\y_{i}\in \RW, i=1,2,\dots, \ell$,
where $\y_0 = {\bf 0}$ and
$$
\widetilde{c}_N = \frac{\widehat{c}_N}{\widehat{c^{\prime}}_N}
= \left( \frac{\pi}{2} \right)^{N/2}
\prod_{j=1}^N \frac{\Gamma (2j)}{\Gamma (j)}.
$$

\vskip 3mm
\noindent(iv)
Consider an ensemble of $N\times N$ complex Hermitian matrices 
$\{ H \}$.
The Gaussian unitary ensemble (GUE) is the ensemble with
the probability density function
$$
\mu^{GUE}(H , \sigma_1^2)
= c_1 \exp \left\{ - \frac{1}{2\sigma_1^2} {\rm Tr}H^2 \right\},
$$
where $\sigma_1^2$ is variance and 
$c_1= 2^{-N/2}(\pi\sigma_1^2)^{-N^2/2}$.
The Gaussian orthogonal ensemble (GOE) is defined as 
the ensemble of $N\times N$ real symmetric matrices $\{ A \}$
with the probability density function
$$
\mu^{GOE}(A , \sigma_2^2)
= c_2 \exp \left\{ - \frac{1}{2\sigma_2^2} {\rm Tr}A^2 \right\},
$$
where $\sigma_2^2$ is variance and 
$c_2= 2^{-N/2}(\pi\sigma_2^2)^{-N(N+1)/4}$.
It is known that the distributions of eigenvalues
$\x = (x_1, x_2, \dots, x_N)$
of these matrix ensembles are given as
$$
g^{GUE}(\x , \sigma_1^2)
= \frac{c_N'}{N!} \sigma_1^{-N^2}
\exp \left\{ - \frac{|\x|^2}{2\sigma_1^2}\right\}
h_N(\x)^2,
$$
and
$$
g^{GOE}(\x , \sigma_2^2)
= \frac{c_N}{N!} \sigma_2^{-N(N+1)/2}
\exp \left\{ - \frac{|\x|^2}{2\sigma_2^2}\right\}
h_N(\x)^,
$$
respectively \cite{Meh91}.
Theorem \ref{thm:VT} and Corollary \ref{cor:VI}
give the relation
$$
g_N^T (0, {\bf 0}, T, \y)= N! g^{GOE}(\y , T),
\quad \y \in \RV,
$$
and
$$
p_N (0, {\bf 0}, t, \y)= N! g^{GUE}(\y , t),
\quad \y \in \RV, \ t>0.
$$
In order to study a Gaussian ensemble of complex Hermitian matrices
intermediate between GUE and GOE, Pandey and Mehta
considered the following probability density functions
with a parameter $\alpha \in [0,1]$
$$
\mu^{PM}(H,\alpha) 
= \int dA \ \mu^{GUE}(H-A, 2\alpha^2 v^2)\mu^{GOE}(A, 2(1-\alpha^2)v^2),
$$
where $v^2 = 1/\{2(1+\alpha^2)\}$ \cite{MP83, PM83}.
They have studied a transition from the GOE to the GUE
observed as $\alpha$ changes from $0$ to $1$.
Let $g^{PM}(\x , \alpha)$ be the probability density function
of eigenvalues in this ensemble of Pandey and Mehta.
We can show the equality \cite{KT02}
$$
\left( \frac{t(2T-t)}{T} \right)^{N/2}
g_N^T \left( 0, {\bf 0}, t, \sqrt{\frac{t(2T-t)}{T}}\x \right)
= N! g^{PM}\left( \x, \sqrt{\frac{T-t}{T}}\right).
$$
It is shown in \cite{KT03} that as a consequence of this
equality, the Harish-Chandra formula for an integral
over the unitary group can be obtained.
Similar argument concerning the relation between the process
$\widehat{{\bf X}}(t)$ and the nonstandard classes of random 
matrices is given in \cite{KTNK*1}.

\vskip 3mm
\noindent (v)
Spohn \cite{Sp87} constructed
nonintersecting Brownian motions on a torus 
and discussed the infinite volume limit to
an infinite system of Dyson-type Brownian motions,
which was also constructed by Dirichlet form technique in 
Osada \cite{Osa96}.
For the present $N$ nonintersecting Brownian motion
${\bf X}(t)$ in a finite time interval $(0, T]$, 
two types of {\it temporally inhomogeneous infinite
particle systems} are obtained by setting
$T=T(N)$ and taking $N \to \infty$. 
If we set $T(N)=2N$ and observe the bulk configuration of
particles at time $t=T(N)+s, -\infty < s \leq 0$,
a spatially homogeneous but temporally inhomogeneous
system is derived in the infinite particle limit,
whose multitime correlation functions have the quaternion
determinantal expressions with sine-kernel.
If we set $T(N)=2N^{1/3}$ and the particle configuration
at time $t=T(N)+s, -\infty < s \leq 0$ around the position
$2N^{2/3}-s^2/4$ is observed, a spatially and temporally
inhomogeneous system is derived in $N \to \infty$, in which
multitime correlation functions are given by the
quaternion determinants with Airy-kernel
\cite{NKT02, KNT*2}.
It is easier to prove the limit theorems for Dyson's Brownian
motion model ${\bf Y}(t)$ corresponding to the above two
kinds of infinite particle limits.
The former limit provides a homogeneous infinite system,
which coincides with the system studied by
Spohn \cite{Sp87}, Osada \cite{Osa96} and
Nagao and Forrester \cite{NF98b}, and the latter does
a temporally homogeneous but spatially inhomogeneous
infinite system, which is related with the process recently
studied by Pr\"ahofer and Spohn \cite{PS02}
and Johansson \cite{Joh02}.
See Nagao \cite{Nagao03} for $N \to \infty$ limit of
the process $\widehat{{\bf X}}(t)$.

\vskip 3mm

\SSC{Proof of Theorems}
\subsection{Proof of Theorem \ref{thm:VT}}

Let $N_{N}(m, \v | \u )$, $\u, \v \in \ZV$,
be the total number of the vicious walks, 
in which the $N$ walkers start from $u_i, i=1,2,\dots,N$,
and arrive at the positions 
$v_i, i=1,2,\dots,N$, at time $m$.
Then the probability that such vicious walks with
fixed end-points are realized 
in all possible random walks started from
the given initial configuration is 
$N_{N}(m, \v |\u )/2^{mN}$, 
which is denoted by $V_{N}(m, \v |\u )$.
We also put
$$
V_{N}(m | \u ) =\sum_{\v \in \ZV} V_{N}(m, \v |\u).
$$
Define a subset of the square lattice ${\bf Z}^{2}$,
$$
{\cal L}_{m}=\{(x,y) \in {\bf Z}^{2}:
x+y=\mbox{even}, \ 0 \leq y \leq m \},
$$
and ${\cal E}_{m}$ be the set of all oriented edges which connect the
nearest-neighbor pairs $((x,y), (x',y'))$ 
of vertices with $y'=y+1$ in ${\cal L}_{m}$.
Then each walk of the $i$-th walker
can be represented as a sequence of successive edges
connecting vertices $(u_{i}, 0)$ and $(v_{i}, m)$ 
on $({\cal L}_{m}, {\cal E}_{m})$, which we call the {\it lattice path}
running from $(u_{i}, 0)$ to $(v_{i}, m)$.
If such lattice paths share a common vertex, they are said to 
intersect. Under the vicious walk condition, 
what we consider is a set of all 
$N$-tuples of {\it nonintersecting lattice paths}.
Let 
$\pi_{0}(\{(u_{i}, 0) \}_{i=1}^{N} \to \{ (v_{i}, m)\}_{i=1}^{N})$
be the set of all $N$-tuples $(\pi_{1}, \dots, \pi_{N})$ of 
nonintersecting lattice
paths, in which $\pi_{i}$ runs from $(u_{i}, 0)$ to $(v_{i}, m)$,
$i=1,2,\dots,N$.
$N_{N}(m,  \v | \u)=
|\pi_{0}(\{(u_{i}, 0)\}_{i=1}^{N} \to \{(v_{i}, m)\}_{i=1}^{N})|$ 
and the Karlin-McGregor formula \cite{KM59_1,KM59_2}
gives
$$
N_{N}(m, \v | \u )
= \det_{1 \leq i, j \leq N}
\Big( |\pi( (u_{j}, 0)\to (v_{i}, m) )|\Big),
$$
where $|A|$ denotes the cardinality of a set $A$ and
$\pi((u_{j},0) \to (v_{i},m))$ the set of lattice paths
from $(u_{j},0)$ to $(v_{i}, m)$.
(Such a determinantal formula is
also known as the Lindstr\"om-Gessel-Viennot formula
in the enumerative combinatorics, see
\cite{Li73,GV85,Stem90}.) 
Since
$
|\pi((u_{j}, 0) \to (v_{i}, m) )|= {m \choose (m+u_{j}-v_{i})/2},
$
we have the binomial determinant 
\begin{equation}
\label{eq:LGV-f}
V_{N}(m , \v | \u )
= 2^{-mN} \det_{1 \leq i, j \leq N} 
\left( {m \choose (m+u_{j}-v_{i})/2} \right).
\end{equation}

For $L >0$ we introduce the following functions:
$$
\phi_L(x)= 2\left[\frac{{L}x}{2}\right], \ x \in {\bf R},
\hbox{ and }
\phi_L(\x)= (\phi_L(x_1), \phi_L(x_2),\dots, \phi_L(x_N)),
\ \x \in {\bf R}^N,
$$
where $[a]$ denotes the largest integer not greater than $a$.
We show the following lemmas.
\vskip 3mm
\begin{lem}
\label{lem:V(t,y|x)}
{\rm (i)}
For $t>0$, $\x \in \ZV$ and $\y\in \RV$
\begin{equation}
\left(\frac{{L}}{2}\right)^{N} 
V_N (\phi_{L^2}(t), \phi_L(\y) | \x)
= c_N' t^{-N^2/2}
h_N\left( \frac{\x}{L}\right)
\exp \left\{ -\frac{|\y|^2}{2t}\right\} h_N(\y)
 \left( 1+ {\cal O}\left(\frac{|\y|}{L}\right) \right),
\label{eq:mto8(1)}
\end{equation}
as $L \to \infty$,
where 
$c^{\prime}_{N}=(2\pi)^{-N/2}
/\prod_{j=1}^{N} \Gamma(j)$.
\par\noindent
{\rm (ii)}
For $t>0$ and $\x \in \ZV$
\begin{equation}
V_N (\phi_{L^2}(t) | \x)
= \frac{1}{\overline{c}_N}
h_N\left(\frac{\x}{L\sqrt{t}}\right) 
\left( 1+ {\cal O}\left(\frac{1}{L}\right) \right),
\label{eq:mto8(2)}
\end{equation}
as $L \to \infty$, where
$\overline{c}_N=\pi^{N/2}
\prod_{j=1}^{N} \{\Gamma(j)/\Gamma(j/2)\}$.
\end{lem}
{\it Proof.} 
It is enough to consider the case that
$\x = 2\u$,  $\phi_L(\y) =2\v$, $\u, \v \in 
{\bf Z}_{<}^{N}$ and $\phi_{L^2}(t) =2\ell$,
$\ell \in {\bf Z}_{+}$.
Then 
$$
N_N (\phi_{L^2}(t) , \phi_L(\y) |\x) 
=
N_N ( 2\ell , 2\v |2\u ) 
= \det_{1 \leq i, j \leq N}
\left( 
{2\ell \choose \ell +u_j-v_i}
\right),
$$
and
\begin{eqnarray}
{2\ell \choose \ell +u_j-v_i}
&=&
\frac{(2\ell)!}{(\ell+u_j-v_i)! (\ell-u_j+v_i)!}
\nonumber
\\
&=&
\frac{(2\ell)!}{(\ell-v_i)! (\ell+v_i)!}A_{ij}(\ell,\v,\u),
\nonumber
\end{eqnarray}
with
$$
A_{ij}(\ell,\v,\u)=
\frac{(\ell+v_i-u_j+1)_{u_j}}{(\ell-v_i+1)_{u_j}},
$$
where 
$(a)_0 \equiv 1$,
$(a)_k= a(a+1) \cdots (a+k-1)$, $k\ge 1$.
Then
\begin{equation}
N_N (\phi_{L^2}(t), \phi_L(\y) |\x) 
= \prod_{i=1}^N
\frac{(2\ell)!}{(\ell-v_i)! (\ell+v_i)!}
\det_{1 \leq i, j \leq N}
\left( 
A_{ij}(\ell,\v,\u)
\right).
\label{eq:Nn}
\end{equation}
The leading term of 
$\det_{1 \leq i, j \leq N}( A_{ij}(\ell,\v,\u))$ 
in $L\to\infty$ is
\begin{eqnarray}
D_1(\v,\u) 
&=& 
\det_{1 \leq i, j \leq N}
\left( \left(\frac{\ell +v_i}{\ell -v_i}\right)^{u_j} \right)
\nonumber
\\
&=&
(-1)^{N(N-1)/2}
\det_{1 \leq i, j \leq N}
\left( \left(\frac{\ell +v_i}{\ell -v_i}\right)^{u_{N-j+1}} \right).
\nonumber
\end{eqnarray}
Let $\xi(\u)=(\xi_{1}(\u), \dots, \xi_{N}(\u))$ 
be a partition specified by the starting point
$2\u$ defined by
\begin{equation}
 \xi_{j}(\u)=u_{N-j+1}-(N-j), \ j =1,2,\dots, N.
\label{eqn:xi}
\end{equation}
Noting that the Vandermonde determinant
$\det_{1\le i,j \le N}(z_i^{N-j})= \prod_{1\le i<j \le N}(z_i -z_j)$,
we have
\begin{eqnarray}
D_1(\v,\u) 
&&= 
(-1)^{N(N-1)/2}
\det_{1 \leq i, j \leq N}
\left( \left(\frac{\ell +v_i}{\ell -v_i}\right)^{N-j} \right)
s_{\xi(\u)}
\left(\frac{\ell+v_1}{\ell-v_1},\dots,\frac{\ell+v_N}{\ell-v_N} \right)
\nonumber
\\
&&=
(-1)^{N(N-1)/2}
\prod_{1\le i<j\le N}
\left( \frac{\ell+v_i}{\ell- v_i}- \frac{\ell+v_j}{\ell-v_j} \right)
s_{\xi(\u)}
\left(\frac{\ell+v_1}{\ell-v_1},\dots,\frac{\ell+v_N}{\ell-v_N} \right)
\nonumber
\\
&&=
\prod_{1\le i<j\le N}
\frac{2\ell(v_j-v_i)}{(\ell-v_i)(\ell-v_j)}
s_{\xi(\u)}
\left(\frac{\ell+v_1}{\ell-v_1},\dots,\frac{\ell+v_N}{\ell-v_N} \right),
\nonumber
\end{eqnarray}
where $s_{\lambda}(z_{1}, \dots, z_{N})$ 
is the Schur function associated to a partition 
$\lambda=(\lambda_{1}, \lambda_{2}, \dots, \lambda_{N})$
defined by 
\begin{equation}
s_{\lambda}(z_{1}, \dots, z_{N}) =
\frac{\det_{1 \leq i, j \leq N} 
\left(z_{i}^{\lambda_{j}+N-j}\right)}
{\det_{1 \leq i, j \leq N} \left(z_{i}^{N-j}\right)}.
\end{equation}
(See Macdonald \cite{Mac95}.)
It is a symmetric polynomial of degree 
$\sum_{i=1}^{N} \lambda_{i}$ in $z_{1}, \dots, z_{N}$ and
it is known that (see p.44 in \cite{Mac95})
\begin{equation}
s_{\lambda}(1,1, \dots, 1)=
\prod_{1 \leq i < j \leq N}
\frac{\lambda_{i}-\lambda_{j}+j-i}{j-i}.
\end{equation}
Then the leading term of $D_1(\v,\u)$ in $L\to\infty$ is
\begin{eqnarray}
D_2(\v,\u) 
&=&
\prod_{1\le i<j\le N}
\frac{2(v_j-v_i)}{\ell}
s_{\xi(\u)}(1,1,\dots,1)
\nonumber
\\
&=&
\ell^{- N(N-1)/2}2^{N(N-1)/2}h_N(\v)h_N(\u)
\prod_{1\le i<j\le N}\frac{1}{j-i}
\nonumber
\\
&=&
h_N\left(\frac{\v}{\ell}\right)h_N(2\u)
\prod_{j=1}^N \frac{1}{\Gamma(j)}.
\label{eq:D2}
\end{eqnarray}
By Stirling's formula we see that
\begin{equation} 
\prod_{i=1}^N
\frac{(2\ell)!}{(\ell-v_i)! (\ell+v_i)!}
= (\ell \pi)^{-N/2} 2^{2N\ell}
\prod_{i=1}^N 
\left( 1- \frac{v_i^2}{\ell^2}\right)^{-\ell-1/2}
\left( \frac{1-v_{i}/\ell}{1+v_{i}/\ell}\right)^{v_i}
\left(1+{\cal O}\left(\frac{1}{\ell}\right)\right).
\label{eq:St}
\end{equation}
From (\ref{eq:Nn}), (\ref{eq:D2}) and (\ref{eq:St})
\begin{eqnarray}
&&V_N(\phi_{L^2}(t), \phi_L(\y) | \x)
= 2^{-2N\ell}N_N(\phi_{L^2}(t), \phi_L(\y) | \x)
\nonumber
\\
&& =
c_N' \left(\frac{2}{\ell}\right)^{N/2}
h_N\left(\frac{\v}{\ell}\right)
h_N\left( 2 \u \right) 
\exp \left\{ - \frac{|\v|^2}{\ell} \right\}
\left(1+{\cal O}\left(\frac{|\v|}{\ell}\right)\right)
\nonumber
\\
&&= 
c_N' \left(\frac{2}{L}\right)^N t^{-N^2/2}
h_N\left(\frac{\x}{L}\right)
\exp \left\{- \frac{|\y|^2}{2t}\right\}
h_N(\y)\left(1+{\cal O}\left(\frac{|\y|}{L}\right)\right).
\nonumber
\end{eqnarray}
Then we obtain (\ref{eq:mto8(1)}).

By (\ref{eq:mto8(1)}) and simple calculation
we have
\begin{equation}
V(\phi_{L^2}(t)|\x)= c_N't^{-N^2/2}
\frac{1}{\Gamma(N+1)}h_N\left(\frac{\x}{L}\right)
\int_{{\bf R}^N} d\y \ e^{- |\y|^{2}/2t}|h_N(\y)|
\left(1+{\cal O}\left(\frac{1}{L}\right)\right),
\end{equation}
as $L\to\infty$.
The last integral is the special case 
($\gamma=1/2$ and $a=1/2t$) of 
\begin{equation}
\int_{{\bf R}^N} d\u \ e^{- a |\u|^{2}}
\prod_{1 \leq i < j \leq N}
|u_{j}-u_{i}|^{2 \gamma} 
= (2 \pi)^{N/2}
(2a)^{-N(\gamma(N-1)+1)/2} 
\prod_{i=1}^{N} \frac{\Gamma(1+i \gamma)}{\Gamma(1+\gamma)}
\label{eq:Sel}
\end{equation}
found in Mehta (eq.(17.6.7) on page 354 in \cite{Meh91}),
whose proof was given in \cite{Mac82}.
Then we have (\ref{eq:mto8(2)}) by elementary calculation.
This completes the proof.
\qed
\vskip 3mm
\begin{lem}
\label{lem:tofn}
Let $t>0$ and $\x, \y \in \RV$. Then
\begin{equation}
\left(\frac{{L}}{2}\right)^{N} 
V_{N} \left( \phi_{L^2}(t) , \phi_L(\y) | \phi_L(\x) \right)
=f_{N}(t, \y |\x )
\left(1+{\cal O}\left(\frac{|\x-\y|}{L}\right)\right),
\end{equation}
as $L\to\infty$.
\end{lem}
{\it Proof.} 
From (\ref{eq:LGV-f})
\begin{eqnarray}
&&\left(\frac{{L}}{2}\right)^{N} 
V_{N} \left( \phi_{L^2}(t) , \phi_L(\y) | \phi_L(\x) \right)
\nonumber
\\
&&= 
2^{-N\phi_{L^2}(t)}\left(\frac{L}{2}\right)^{N} 
\det_{1 \leq i, j \leq N}
\left( \phi_{L^2}(t) \choose (\phi_{L^2}(t)+\phi_L(x_j)-\phi_L(y_i))/2
\right)
\nonumber\\
&& =
\det_{1 \leq i, j \leq N}
\left( 
2^{-\phi_{L^2}(t) -1} L
{\phi_{L^2}(t) \choose 
(\phi_{L^2}(t)+\phi_L(x_j)-\phi_L(y_i))/2}\right).
\nonumber
\end{eqnarray}
Application of Stirling's formula yields the lemma.
\qed
\vskip 3mm

By Donsker's theorem (see, for instance, Billingsley \cite{Bi68}) 
we see that ${\cal N}_{N}(t, \x )$ 
is the probability that $N$ Brownian motions do not collide until time $t$.
We have the following asymptotic behaviours of the function
${\cal N}_{N}(t, \x)$ as $|\x|/\sqrt{t} \to 0$.

\begin{lem}
\label{lem:N(t,x)}
Let $t > 0$ and $\x \in \RV$. Then
\begin{equation}
{\cal N}_{N}(t, \x)
=\frac{1}{\overline{c}_N}
h_N \left( \frac{\x}{\sqrt{t}}\right)
\left( 1+ {\cal O}\left(\frac{|\x|}{\sqrt{t}}\right) \right),
\quad \frac{|\x|}{\sqrt{t}} \to 0,
\label{tto8}
\end{equation}
where $\overline{c}_{N}=\pi^{N/2} \prod_{j=1}^{N}
\{\Gamma(j)/\Gamma(j/2)\}$.
\end{lem}
{\it Proof.} 
First note that
$$
f_N (t, \y | \x)= (2 \pi t)^{-N/2} 
\exp\left\{ -\frac{1}{2t}\sum_{i=1}^N (x_{i}^2+y_{i}^2) \right\}
\det_{1 \leq i, j \leq n}
\left( e^{x_{j} y_{i}/t} \right). \nonumber
$$
We rewrite the determinant as
\begin{eqnarray}
\det_{1 \leq i, j \leq N}
\left( e^{x_{j} y_{i}/t} \right)
&=& \frac{\det_{1 \leq i, j \leq N} \left(
(e^{x_{i}/t})^{y_{N-j+1}} \right)}
{\det_{1 \leq i, j \leq N} \left(
(e^{x_{i}/t})^{N-j} \right)}
\times 
\det_{1 \leq i, j \leq N} \left((e^{x_{i}/t})^{N-j}\right)
\nonumber
\\
&=& s_{\xi(\y)}
\left(e^{x_{1}/t}, e^{x_{2}/t},\dots,e^{x_{N}/t}\right)
\prod_{1 \leq i < j \leq N}(e^{x_{j}/t}-e^{x_{i}/t}),
\nonumber
\end{eqnarray}
where $\xi_i(\y)= y_{N-i+1}-(N-i)$, $i=1,2,\dots,N$.
Using it
\begin{eqnarray}
f_{N}(t, \y|\x)
&=& (2 \pi t)^{-N/2} 
s_{\xi(\y)}\left(e^{x_{1}/t}, e^{x_{2}/t}, \dots,
e^{x_{N}/t}\right)
\nonumber
\\
&& \ \times \exp\left\{- \frac{|\x|^2 + |\y|^2}{2t} \right\}
\prod_{1 \leq i < j \leq N}
(e^{x_{j}/t}-e^{x_{i}/t}).
\end{eqnarray}
Since 
$$
\lim_{\frac{|\x |}{{t}} \to 0} 
s_{\xi(\y)}(e^{x_{1}/t}, \dots, e^{x_{N}/t})
= s_{\xi(\y)}(1,1, \dots, 1) 
= h_N(\y)\prod_{j=1}^N \frac{1}{\Gamma (j)},
$$
and
$$ 
\prod_{1 \leq i < j \leq N}
(e^{x_{j}/t}-e^{x_{i}/t})
= h_N \left(\frac{\x}{t}\right)
\left( 1+ {\cal O}\left(\frac{|\x|}{{t}}\right) \right),
\quad \frac{|\x|}{{t}} \to 0,
$$
the function is asymptotically
\begin{eqnarray}
{\cal N}_{N}(t, \x ) 
&=& (2 \pi)^{-N/2} t^{-N(N+1)/4}
h_N \left(\frac{\x}{\sqrt{t}}\right)
\prod_{j=1}^N \frac{1}{\Gamma (j)}
\nonumber
\\
&&\times\int_{\RV} d\y 
h_N(\y)\exp\left\{ - \frac{|\y |^2 }{2t} \right\}
\left( 1+ {\cal O}\left(\frac{|\x|}{\sqrt{t}}\right) \right),
\quad \frac{|\x|}{\sqrt{t}}\to 0.
\nonumber
\end{eqnarray}
By (\ref{eq:Sel}) we have (\ref{tto8}).
\qed
\vskip 3mm
\begin{lem}
\label{lem:tight}
For $\z \in \ZV$ and $T>0$
$\{ \mu_{L,T}^{\z}, L \ge 1\}$ is tight.
\end{lem}
{\it Proof.} 
By the Kolmogorov's tightness criterion
it is enough to prove that for any $\varepsilon>0$
\begin{equation}
\lim_{\delta\to 0}\sup_{L\ge 1}
\mu_{L,T}^{\z}
\left( \max_{0\le u, v \le T, |u-v|<\delta}
|w(u)- w(v)| \ge \varepsilon \right)
=0.
\label{Kolmogorov}
\end{equation}
(See, for example, Billingsley \cite{Bi68}.)
Since $\mu_{L,T}^{\z}$ is the probability measure 
of a linearly interpolated random process,
(\ref{Kolmogorov}) is derived from the following estimates:
as $\delta \to 0$,
\begin{equation}
\limsup_{L\to \infty}
\mu_{L,T}^{\z}
\left( \max_{0\le u \le \delta}
|w(u)- w(0)| \ge \varepsilon/2 \right)
= {o (\delta)},
\label{tihght1}
\end{equation}
\begin{equation}
\limsup_{L\to \infty}
\mu_{L,T}^{\z}
\left(  \max_{0\le u \le \delta}
|w(t+u)- w(t)| \ge \varepsilon/2 \right)
= {o (\delta)},
\quad t\in [\delta/2, T-\delta].
\label{tihght2}
\end{equation}
Under the nonintersecting condition,
for any $i=1,2,\dots,N$
\begin{equation}
|w_i(u)-w_i(0)|\le 
|w_N(0)-w_1(0)|+ (w_1(u)-w_1(0))_{-}+ (w_N(u)-w_N(0))_{+},
\nonumber
\end{equation}
where $a_+ = \max \{a,0\}$ and $a_- = \max \{-a,0\}$.
Then the set
$$
\left\{ \max_{0\le u \le \delta} |w(u)- w(0)| \ge \varepsilon/2 \right\}
$$
is included in the set
$$
\left\{ \max_{0\le u \le \delta} (w_1(u)- w_1(0))_{-} 
\ge \frac{\varepsilon}{4}
-\frac{z_N-z_1}{2{L}} \right\}
\cup
\left\{ \max_{0\le u \le \delta} (w_N(u)- w_N(0))_{+} 
\ge \frac{\varepsilon}{4}
-\frac{z_N-z_1}{2{L}} \right\}.
$$
Noting that $(w_1(u)- w_1(0))_{-}$ and $(w_N(u)- w_N(0))_{+}$
are nonnegative submartingales, 
we can apply Doob's theorem 
(see, for instance, Revus and Yor \cite{RY98})
to obtain
\begin{eqnarray}
&&\mu_{L,T}^{\z}\left( \max_{0\le u \le \delta}
|w(u)- w(0)| \ge \varepsilon/2 \right)
\nonumber
\\
&&\le \left(\frac{8}{\varepsilon}\right)^p
E_{L,T}^{\z}\Big( |w_1(\delta)- w_1(0)|^p +|w_N(\delta)- w_N(0)|^p\Big)
\end{eqnarray}
for any $p>1$ and $L > 4(z_N-z_1)/\varepsilon $, where 
$E_{L,T}^{\z}$ represents the expectation with respect to
the probability measure $\mu_{L,T}^{\z}$.

From Lemmas \ref{lem:V(t,y|x)}, \ref{lem:tofn}
and \ref{lem:N(t,x)}
\begin{eqnarray*}
&&\limsup_{L\to\infty}E_{L,T}^{\z}
\Big( |w_1(\delta)- w_1(0)|^p +|w_N(\delta)- w_N(0)|^p \Big)
\\
&&\le C_1 \limsup_{L\to\infty}
\int_{\RV}d\y\
\left(\frac{{L}}{2}\right)^N
\frac{V_N(\phi_{L^2}(\delta),\phi_L(\y)|\z)
V_N(\phi_{L^2}(T-\delta)|\phi_L(\y) )}{V_N(\phi_{L^2}(T)|\z)}
( y_1^p + y_N^p)
\\
&&\le C_2 c_N \delta^{-N^2/2}
\int_{\RV}d\y \ 
\exp\left\{-\frac{|\y|^2}{2\delta} \right\}h_N(\y)^2( y_1^p + y_N^p)
\\
&&\le C_3 \delta^{p/2}
\int_{\RV}d\x \ \exp\left\{-\frac{|\x|^2}{2} \right\}h_N(\x)^2( x_1^p + x_N^p)
\\
&& ={\cal O}(\delta^{p/2}).
\end{eqnarray*}
Taking $p> 2$, we obtain (\ref{tihght1}).

Fix $t\in [\delta/2, T-\delta]$.
By the Markov property
\begin{eqnarray}
&&\mu_{L,T}^{\z}\left( \max_{0\le u \le \delta}
|w(t+u)- w(t)| \ge \varepsilon/2 \right)
\nonumber
\\
&&=Q^{\z}_{L^2T}\left(\max_{0\le u \le \delta}
\left|\frac{S(L^2(t+u))-S(L^2t)}{L}\right| \ge \varepsilon/2 \right)
\nonumber
\\
&&\le \frac{1}{P^{\z} ( \Lambda_{L^2T})}
E^{\z}\left(\Lambda_{L^2t},
P^{S(L^2t)}\left(\max_{0\le u \le \delta}
\left|\frac{S(L^2u)-S(0)}{L}\right| \ge \varepsilon/2
\right)\right),
\end{eqnarray}
where 
$\Lambda_m = \{S_j^1 < S_j^2 < \cdots < S_j^N, \; 0\le j \le m \}$
and $E^{\z}$ represents the expectation with respect to 
the probability measure $P^{\z}$.
By Doob's inequality for any $\x\in {\bf Z}^N$ and $p>1$
\begin{equation}
P^{\x}\left(\max_{0\le u \le \delta}
\left|\frac{S(L^2u)-S(0)}{{L}}\right| \ge \varepsilon/2\right)
\le \left(\frac{2}{\varepsilon}\right)^p 
E^{\x}\left( \left|\frac{S(L^2\delta)-S(0)}{{L}}\right|^p\right)
\le C_4 \delta^{p/2}.
\end{equation}
From Lemma \ref{lem:V(t,y|x)} {\rm (ii)}
\begin{eqnarray}
&&\limsup_{L\to\infty}
\mu_{L,T}^{\z}\left( \max_{0\le u \le \delta}
|w(t+u)- w(t)| \ge \varepsilon/2 \right) \nonumber\\
&&\quad \le C_4 \delta^{p/2}\limsup_{L\to\infty}
\frac{P^{\z}(\Lambda_{L^2\delta/2})}{P^{\z}(\Lambda_{L^2T})}
={\cal O}(\delta^{p/2- N(N-1)/4}).
\end{eqnarray}
Taking $p>N(N-1)/2 +2$, we obtain (\ref{tihght2}).
This completes the proof.
\qed
\vskip 3mm
\begin{lem}
\label{lem:FDD}
Let $\z\in\ZV$, $0=t_0<t_1< \cdots <t_{k-1}<t_k=T$ and
$\theta=(\theta_1, \dots, \theta_k)\in {\bf R}^{Nk}$.
Then
\begin{eqnarray}
&&\lim_{L\to\infty}
E_{L,T}^{\z}
\left(\exp\left\{\sqrt{-1}\sum_{j=1}^k \theta_j\cdot w(t_j)\right\}\right)
\nonumber
\\
&&= \int_{(\RV)^k} d\y_1d\y_2 \cdots d\y_k 
\prod_{j=1}^k g^T_N(t_{j-1},\y_{j-1},t_j,\y_j)
\exp\left\{\sqrt{-1}\sum_{j=1}^k \theta_j\cdot \y_j \right\},
\end{eqnarray}
\end{lem}
where $\y_0 = {\bf 0}$.

{\it Proof.} 
By Lemmas \ref{lem:tofn} and \ref{lem:V(t,y|x)}
\begin{eqnarray}
&&\lim_{L\to\infty}
E_{L,T}^{\z}
\left(\exp\left\{\sqrt{-1}\sum_{j=1}^k \theta_j\cdot w(t_j)\right\}\right)
\nonumber
\\
&&=\lim_{L\to\infty}
\frac{1}{P^{\z}(\Lambda_{L^2T})}
\sum_{\x_1\in\ZV} \cdots \sum_{\x_k\in\ZV}
E^{\z}\Bigg[ \Lambda_{L^2t_1}, S(\phi_{L^2}(t_1))=\x_1,
\nonumber
\\
&&\qquad\times
E^{\x_1}\Bigg[ \Lambda_{L^2(t_2-t_1)}, S(\phi_{L^2}(t_2-t_1))=\x_2,
\cdots
\nonumber
\\
&&\qquad \times
\left.\left.
E^{\x_{k-1}}\left[ \Lambda_{L^2(t_{k}-t_{k-1})}, 
S(\phi_{L^2}(t_{k}-t_{k-1})=\x_{k},
\exp\left\{\sqrt{-1}\sum_{j=1}^k \theta_j\cdot \frac{\x_j}{{L}}\right\}
\right]\right]\right]
\nonumber
\\
&&=\lim_{L\to\infty}
\frac{1}{P^{\z}(\Lambda_{L^2T})}
\int_{(\RV)^k} d\y_1d\y_2 \cdots d\y_k \
\left(\frac{{L}}{2}\right)^{Nk}
V_N (\phi_{L^2}(t_1), \phi_L(\y_1)|\z)
\nonumber
\\
&&\qquad \times
V_N (\phi_{L^2}(t_2-t_1), \phi_L(\y_2)|\phi_L(\y_1))\cdots
\nonumber
\\
&&\qquad \times
V_N (\phi_{L^2}(t_k-t_{k-1}), \phi_L(\y_k)|\phi_L(\y_{k-1}))
\exp\left\{\sqrt{-1}\sum_{j=1}^k \theta_j\cdot \y_j \right\}
\nonumber
\\
&&= \int_{(\RV)^k} d\y_1d\y_2 \cdots d\y_k 
\prod_{j=1}^k g^T_N(t_{j-1},\y_{j-1},t_j,\y_j)
\exp\left\{\sqrt{-1}\sum_{j=1}^k \theta_j\cdot \y_j \right\}.
\nonumber
\end{eqnarray}
This completes the proof.
\qed
\vskip 3mm
{\it Proof of Theorem \ref{thm:VT}.}
From Lemmas \ref{lem:tight} and \ref{lem:FDD}
we see that 
$\mu_{L,T}^{\z}(\cdot)$ converges weakly to the law of 
the time inhomogeneous diffusion process ${\bf X}(t)$
with transition density $g^T_N(s,\x,t,\y)$.
Noting that 
${\cal N}_N(T-t,\x)$ is a solution of the heat equation,
we see that
$g_N^T(s,\x, t, \y )$ satisfies the following backward equation:
\begin{equation}
\frac{\partial}{\partial t}g_N^T(s,\x, t, \y )
= \frac{1}{2}\triangle_{\x}g_N^T(s,\x, t, \y )
+ b^T(t,\x)\cdot \nabla_{\x}g_N^T(s,\x, t, \y ).
\end{equation}
Then the process ${\bf X}(t)$ solves (\ref{eq:SDE1}).
This completes the proof of Theorem \ref{thm:VT}.
\vskip 3mm
{\it Proof of Corollary \ref{cor:VI}.}
By simple observation we see that
the estimates concerning the probability $\mu^{\z}_{L,T}$ 
in the lemmas are uniform with respects to $T$.
Then (i) is obtained from the properties
\begin{eqnarray*}
&&\lim_{T\to \infty}g_N^T(0,{\bf 0}, t, \y )
= p_N (0,{\bf 0}, t, \y ),
\\
&&\lim_{T\to \infty}g_N^T(s,\x, t, \y )
=p_N (s,\x, t, \y ),
\end{eqnarray*}
which are derived from Lemma \ref{lem:N(t,x)}.
Noting that 
$h_N(\x)$ is a harmonic function,
we see that
$p_N (s,\x, t, \y )$ satisfies the following backward equation:
\begin{equation}
\frac{\partial}{\partial t}p_N(s,\x, t, \y )
= \frac{1}{2}\triangle_{\x}p_N(s,\x, t, \y )
+ \nabla_{\x}\ln h_N(\x)\cdot \nabla_{\x}p_N(s,\x, t, \y ).
\end{equation}
Then the process ${\bf Y}(t)$ solves (\ref{eq:SDE2}).
This completes the proof of Corollary \ref{cor:VI}.
\qed

\subsection{Proof of Theorem \ref{thm:WT}}

Let $\widehat{N}_{N}(m, \v | \u )$, $\u, \v \in \ZW$, 
be the total number of the vicious walks with wall restriction, 
in which the $N$ walkers 
start form the positions $u_i, i=1,2,\dots,N$, and
arrive at the positions $v_i, i=1,2,\dots,N$, at time $m$.
Then the probability that such vicious walks with
fixed end-points are realized 
in all possible random walks started from
the given initial configuration is 
$\widehat{N}_{N}(m, \v |\u )/2^{mN}$, 
which is denoted by $\widehat{V}_{N}(m, \v |\u )$.
We also put
$$
\widehat{V}_{N}(m | \u ) =\sum_{\v \in \ZW} \widehat{V}_{N}(m, \v |\u).
$$
By the Karlin-McGregor (Lindstr\"om-Gessel-Viennot) formula, 
we have \cite{KGV00}
\begin{equation}
\label{eq:LGV_WW}
\widehat{V}_{N}(m,\v | \u)
= 2^{-mN}\det_{1 \leq i, j \leq N} 
\left(
{ m \choose (m+u_{j}-v_{i})/2 }
- { m \choose (m+u_{j}+v_{i})/2+1 }
\right).
\end{equation}
Let 
\begin{equation}
sp_{\lambda}(z_{1}, \dots, z_{N})=
\frac{\det(z_{i}^{\lambda_{j}+N-j+1}-z_{i}^{-(\lambda_{j}+N-j+1)})}
{\det(z_{i}^{N-j+1}-z_{i}^{-(N-j+1)})},
\end{equation}
for a partition $\lambda=(\lambda_{1}, \dots, \lambda_{N})$.
Remark that $sp_{\lambda}(z_{1}, \dots, z_{N})$ is 
the character of the irreducible representation 
corresponding to a partition $\lambda$
of the symplectic Lie algebra 
(see, for example, Lectures 6 and 24 in 
Fulton and Harris \cite{FH91}). 
By using the function $sp_{\lambda}$ instead of the
Schur function, we can show the following lemma
by a similar way to the proof of Lemma \ref{lem:V(t,y|x)}.

\vskip 3mm
\begin{lem}
\label{lem:hV(t,y|x)}
{\rm (i)}
For $t>0$, $\x \in \ZW$ and $\y\in \RW$
\begin{equation}
\left(\frac{{L}}{2}\right)^{N} 
\widehat{V}_N (\phi_{L^2}(t), \phi_L(\y) | \x)
= \widehat{c^{\prime}}_N t^{-N(2N+1)/2}
\h \left(\frac{\x}{L} \right)
\exp \left\{ -\frac{|\y|^2}{2t}\right\} \h (\y)
 \left( 1+ {\cal O}\left(\frac{|\y|}{L}\right) \right),
\label{eq:hVmto8(1)}
\end{equation}
as $L \to \infty$, where
$\widehat{c^{\prime}}_N=(2/\pi)^{N/2}
/ \prod_{j=1}^{N} \Gamma(2j)$.
\par\noindent
{\rm (ii)}
For $t>0$ and $\x \in \ZW$
\begin{equation}
\widehat{V}_N (\phi_{L^2}(t) | \x)
= \frac{1}{\widetilde{c}_N}
\h\left(\frac{\x}{L\sqrt{t}}\right) 
\left( 1+ {\cal O}\left(\frac{1}{L}\right) \right),
\label{eq:hVmto8(2)}
\end{equation}
as $L \to \infty$, where
$\widetilde{c}_N=(\pi/2)^{N/2}
\prod_{j=1}^{N} \{\Gamma(2j)/\Gamma(j)\}$.
\end{lem}
{\it Proof.} 
Again we consider the case that
$\x = 2\u$,  $\phi_L(\y) =2\v$,
$\u, \v \in {\bf N}_{<}^{N}$ and $\phi_{L^2}(t) =2\ell$,
$\ell \in {\bf Z}_{+}$.
By the equation (\ref{eq:LGV_WW})
\begin{equation}
\label{eq:LGV_WW2}
\widehat{V}_{N}(\phi_{L^2}(t) , \phi_L (\y) | \x )
= 2^{-2N\ell}
\prod_{i=1}^N \frac{(2\ell)!}{(\ell-v_i)! (\ell+v_i)!}
\det_{1 \leq i, j \leq N} 
\left(
\widehat{A}_{ij}(\ell,\v,\u)
\right),
\end{equation}
with
$$
\widehat{A}_{ij}(\ell,\v,\u)
=
\frac{(\ell+v_i-u_j+1)_{u_j}}{(\ell-v_i+1)_{u_j}}
-
\frac{(\ell-v_i-u_j+1)_{u_j+1}}{(\ell+v_i+1)_{u_j+1}}.
$$
Then the leading term of 
$\det_{1 \leq i, j \leq N}( \widehat{A}_{ij}(\ell,\v,\u))$ 
as $L\to\infty$ is
\begin{eqnarray}
&&\widehat{D}_1(\v,\u) 
=\det_{1 \leq i, j \leq N}
\left( \left(\frac{\ell +v_i}{\ell -v_i}\right)^{u_j} 
- \left(\frac{\ell -v_i}{\ell +v_i}\right)^{u_j}
\right)
\nonumber
\\
&&\qquad\qquad=
(-1)^{N(N-1)/2}
\det_{1 \leq i, j \leq N}
\left( \left(\frac{\ell +v_i}{\ell -v_i}\right)^{u_{N-j+1}} 
- \left(\frac{\ell -v_i}{\ell +v_i}\right)^{u_{N-j+1}} 
\right)
\nonumber
\\
&&\qquad\qquad=
\det_{1 \leq i, j \leq N}
\left( 
\left( \frac{\ell -v_i}{\ell +v_i}\right)^{N-j+1}
-\left(\frac{\ell +v_i}{\ell -v_i}\right)^{N-j+1} 
\right)
sp_{\widehat{\xi}(\u)}
\left( \frac{\ell+v_1}{\ell-v_1},\dots,
\frac{\ell+v_N}{\ell-v_N} \right),
\nonumber
\end{eqnarray}
where $\widehat{\xi}(\u)=(\widehat{\xi}_{1}(\u), \dots,
\widehat{\xi}_{N}(\u))$ with
$\widehat{\xi}_{j}(\u)=u_{N-j+1}-(N-j+1),
j=1,2, \dots, N$.
Note that
$$
\det_{1 \leq i, j \leq N}
\left( 
z_i^{N-j+1} - z_i^{-(N-j+1)}
\right)
= \prod_{j=1}^N \left(z_j - \frac{1}{z_j}\right)
\prod_{1\le i<j\le N}\left(z_j-z_i\right)\left(\frac{1}{z_i z_j}-1\right).
$$
Then by simple calculation we have
\begin{eqnarray}
&&\det_{1 \leq i, j \leq N}
\left( 
\left(\frac{\ell -v_i}{\ell +v_i}\right)^{N-j+1}
- \left(\frac{\ell +v_i}{\ell -v_i}\right)^{N-j+1} 
\right)
\nonumber
\\
&&\qquad=
\prod_{j=1}^{N}
\frac{4\ell v_j}{\ell^2 - v_j^2}
\prod_{1 \leq i < j \leq N} 
\frac{4\ell^2 (v_j^2 -v_i^2)}{(\ell^2 -v_i^2)(\ell^2 -v_j^2)}.
\end{eqnarray}
It is known that
\begin{equation}
\label{eq:sp(11)}
sp_{\lambda}(1, \dots, 1)=
\prod_{1 \leq i < j \leq N} 
\frac{\ell_{j}^2-\ell_{i}^{2}}{m_{j}^{2}-m_{i}^{2} }
\prod_{j=1}^{N} 
\frac{\ell_{j}}{m_{j}}
\end{equation}
with
$\ell_{j}=\lambda_{j}+N-j+1$, $m_{j}=N-j+1$
\cite{FH91}.
Then 
\begin{equation}
sp_{\widehat{\xi}(\u)}(1, \dots, 1)
=
\prod_{1 \leq i < j \leq N} 
\frac{u_j^2 - u_i^2}{j^2-i^2}
\prod_{j=1}^{N} 
\frac{u_j}{j}
= \h (\u)\prod_{j=1}^N \frac{1}{\Gamma (2j+1)}.
\end{equation}
Then the leading term of $\widehat{D}_1(\v,\u)$ in $L\to\infty$ is
\begin{eqnarray}
&&\widehat{D}_2(\v,\u)=
\h (\u)\prod_{j=1}^N \frac{1}{\Gamma (2j)}
\prod_{j=1}^{N}
\frac{4 v_j}{\ell}
\prod_{1 \leq i < j \leq N} 
\frac{4 (v_j^2 -v_i^2)}{\ell^2}
\nonumber
\\
&&=\prod_{j=1}^N \frac{2}{\Gamma (2j)}
\h \left(\frac{\u}{\ell}\right) \h (2 \v).
\label{eq:hD_2}
\end{eqnarray}
From (\ref{eq:St}),(\ref{eq:LGV_WW2}) and (\ref{eq:hD_2}) we have
\begin{eqnarray}
&&\widehat{V}_N(\phi_{L^2}(t), \phi_L(\y) | \x)
\nonumber
\\
&& =
\prod_{j=1}^N \frac{1}{\Gamma (2j)}
\left(\frac{4}{\ell\pi}\right)^{N/2}
\h\left(\frac{\v}{\ell}\right)\h(2\u)
\exp \left\{ - \frac{|\v|^2}{\ell} \right\}
\left(1+{\cal O}\left(\frac{|\v|}{\ell}\right)\right)
\nonumber
\\
&&= \left(\frac{2}{L}\right)^N
\widehat{c^{\prime}}_N t^{-N(2N+1)/2}
\h \left( \frac{\x}{L}\right)\h(\y)
\exp \left\{ - \frac{|\y|^2}{2t} \right\}
\left(1+{\cal O}\left(\frac{|\y|}{L}\right)\right).
\nonumber
\end{eqnarray}
Then we obtain (\ref{eq:hVmto8(1)}).

By (\ref{eq:hVmto8(1)}) and simple calculation
we have
\begin{eqnarray}
&&\widehat{V}(\phi_{L^2}(t)|\x)
= \frac{\widehat{c^{\prime}}_N t^{-N(2N+1)/2}}{\Gamma(N+1)}
\h\left(\frac{\x}{L}\right)
\int_{{\bf R}^N} d\y \ e^{- |\y|^{2}/2t}|\h(\y)|
\left(1+{\cal O}\left(\frac{1}{L}\right)\right)
\nonumber
\\
&&\qquad
=\frac{\widehat{c^{\prime}}_N}{\Gamma(N+1)}
\h\left(\frac{\x}{L\sqrt{t}}\right)
\int_{{\bf R}^N} d\z \ e^{- |\z|^{2}/2}|\h(\z)|
\left(1+{\cal O}\left(\frac{1}{L}\right)\right),
\nonumber
\end{eqnarray}
as $L\to\infty$.
The last integral is the special case 
($\gamma=1/2$ and $a=1$) of 
\begin{eqnarray}
\label{eq:Mehta2}
&& \int_{{\bf R}^N} d \u e^{- |\u|^{2}/2}
\prod_{1 \leq i < j \leq N}
|u_{j}^{2}-u_{i}^{2}|^{2 \gamma}
\prod_{j=1}^{N} |u_{j}|^{2 a-1} 
\nonumber
\\
&& \qquad\qquad\ = 2^{aN + \gamma N (N-1)}
\prod_{j=1}^{N} \frac{\Gamma(1+j \gamma) 
\Gamma(a+\gamma(j-1))}{\Gamma(1+\gamma)},
\end{eqnarray}
((17.6.6) on page 354 in \cite{Meh91}).
Then we obtain (\ref{eq:hVmto8(2)})
by elementary calculation.
\qed
\vskip 3mm

Following the same calculation as was done in the proof of
Lemma \ref{lem:tofn}, we have the following lemma.
\vskip 3mm

\begin{lem}
\label{lem:tohfn}
Let $t>0$ and $\x, \y \in \RW$.
Then
\begin{equation}
\left(\frac{{L}}{2}\right)^{N} 
\widehat{V}_{N} \left( \phi_{L^2}(t) , \phi_L(\y) | \phi_L(\x) \right)
=\f (t, \y |\x )
\left(1+{\cal O}\left(\frac{|\x|+|\y|}{L}\right)\right),
\end{equation}
as $L\to\infty$.
\end{lem}
\vskip 3mm

We rewrite $\f (t, \y |\x )$ as
\begin{eqnarray}
\f (t, \y |\x )
&=& (2 \pi t)^{-N/2} 
sp_{\widehat{\xi}(\y)}\left(e^{x_{1}/t}, 
\dots, e^{x_{N}/t}\right)
\exp \left\{ - \frac{|\x|^2 + |\y|^2}{2t} \right\}
\nonumber
\\
&\times& 
\prod_{j=1}^{N}(e^{x_{j}/t}-e^{-x_{j}/t}) 
\left\{ \prod_{1 \leq i < j \leq N}
(e^{x_{j}/t}-e^{x_{i}/t})(e^{(x_{i}+x_{j})/t}-1) \right\}
\left\{ \prod_{j=1}^{N} e^{x_{j}/t} \right\}^{-N+1}.
\nonumber
\end{eqnarray}
Then we can obtain the following lemma by a similar way to prove
Lemma \ref{lem:N(t,x)} by virtue of the equations
(\ref{eq:sp(11)}) and (\ref{eq:Mehta2}). 
\vskip 3mm
\begin{lem}
\label{lem:hN(t,x)}
Let $t > 0$ and $\x \in \RW$. Then
\begin{equation}
\hN (t, \x)
=\frac{1}{\widetilde{c}_N}
\h \left( \frac{\x}{\sqrt{t}}\right) 
\left( 1+ {\cal O}\left(\frac{|\x|}{\sqrt{t}} \right) \right),
\quad \frac{|\x|}{\sqrt{t}} \to 0,
\label{Wtto8}
\end{equation}
where
$\widetilde{c}_{N}=(\pi/2)^{N/2}
\prod_{j=1}^{N} \{\Gamma(2j)/\Gamma(j)\}$.
\end{lem}
\vskip 3mm
From the above lemmas and the same argument as in the previous subsection
we can obtain Theorem \ref{thm:WT} and Corollary \ref{cor:WI}.

\vskip 3mm

\SSC{Proof of Lemmas \ref{lem:Nor_V} and \ref{lem:Nor_W}}

We use the following identity, 
which is shown in de Bruijn \cite{deBr55}.
(See also Appendix in \cite{KT02}.)
Lemmas \ref{lem:Nor_V} and \ref{lem:Nor_W}
are easy consequences of this result
as shown below.

\begin{lem}
\label{lem:IOka}
Let $z$ be a square integrable piecewise continuous function on ${\bf R}^2$.
Then
\begin{equation}
\int_{{\RV}} d\y \det_{1 \le i,j \le N}(z(x_i, y_j))
= {\rm Pf}_{1 \leq i < j \leq \hat{N}}( F_{ij}(\x) ),
\end{equation}
where
$$
\hat{N}= \left\{ \begin{array}{cl}
N, & \qquad \mbox{if $N$ is even} \\
N+1, & \qquad \mbox{if $N$ is odd},\\
\end{array} \right.
$$
\begin{eqnarray}
&&I_z(x_{i}) = \int_{-\infty}^{\infty} z(x_{i},y )dy,
\nonumber
\\
&&I_z (x_{i}, x_{j}) = 
\int_{- \infty < y_{1} < y_{2} < \infty} 
\det
\left( \matrix{
z(x_{i},y_{1}) & z(x_{i},y_{2}) \cr
z(x_{j},y_{1}) & z(x_{j},y_{2}) } \right)
dy_{1} dy_{2},
\nonumber
\end{eqnarray}
and 
$$
F_{ij}(\x) =
\left\{ \begin{array}{cl}
I_z(x_{i}, x_{j}),
& \qquad \mbox{if} \quad 1 \leq i < j \leq N, \\
-I_z(x_{i}, x_{j}), 
& \qquad \mbox{if} \quad 1 \leq j < i \leq N, \\
I_z(x_i), & \qquad
\mbox{if} \quad 1 \leq i \leq N, j=N+1, \\
-I_z(x_j), & \qquad
\mbox{if} \quad i=N+1, 1 \leq j \leq N, \\
0, & \qquad
\mbox{if} \quad 1 \leq i=j \leq N+1. \\
\end{array} \right.
$$
\end{lem}

\vskip 3mm

{\it Proof of Lemma \ref{lem:Nor_V}. }
From the above lemma and integration by substitution,
it is enough to show
\begin{equation}
\label{eq:z_1}
I_{z_1}(x_i)=1,
\quad 
I_{z_1}(x_i,x_j)= \Psi \left( \frac{x_j-x_i}{\sqrt{2}}\right),
\end{equation}
for
$z_1(x,y)=e^{-(x -y)^2}/\sqrt{\pi}.$
The first equation in (\ref{eq:z_1}) is trivial.
Let $z_{0}(x)=e^{-x^2}/\sqrt{\pi}$. Then
\begin{eqnarray}
I_{z_1}(x_i,x_j)
&=& \int_{-\infty}^{\infty} dy_1 \int_{x_i-x_j}^{x_j-x_i}dy_2
\; z_0(y_1)z_0(y_1+y_2)
\nonumber
\\
&=& \frac{1}{\pi}\int_{x_i-x_j}^{x_j-x_i}dy_2 \; e^{-y_2^2/2}
\int_{-\infty}^{\infty}  dy_1 \; e^{-2(y_1 + y_2/2)^2}
\nonumber
\\
&=& \frac{1}{\sqrt{2\pi}}\int_{x_i-x_j}^{x_j-x_i}dy_2 \; e^{-y_2^2/2}
=\Psi \left( \frac{x_j-x_i}{\sqrt{2}}\right).
\nonumber
\end{eqnarray}
This completes the proof.
\qed

{\it Proof of Lemma \ref{lem:Nor_W}. }
It is enough to show
\begin{equation}
\label{eq:z_2}
I_{z_2}(x_i)=\Psi \left(x_i \right),
\quad 
I_{z_2}(x_i,x_j)=\widehat{\Psi} \left(x_i, x_j \right),
\end{equation}
for
$z_2(x,y)= \left\{
\left( e^{-(x -y)^2}-e^{-(x +y)^2}\right)/\sqrt{\pi} \right\}
{\bf 1}\{x\ge 0, y\ge 0 \}.
$
The first equation in (\ref{eq:z_2}) is obtained easily.
To show the second equation we put
$$
G((a_1,a_2], (b_1,b_2])
=\int_{a_1}^{a_2}dy_1 \int_{b_1}^{b_2} dy_2
\; z_0(y_1)z_0(y_1 -y_2), 
$$
for $a_1, a_2, b_1, b_2 \in {\bf R}\cup \{-\infty, \infty \}$.
Then 
\begin{eqnarray}
I_{z_2}(x_i,x_j)
&=& G((-x_i,\infty), (-\infty, x_j-x_i])-G((-x_i,\infty), (-\infty, -x_j-x_i])
\nonumber
\\
&-& G((x_i,\infty), (-\infty, x_j+x_i])+G((x_i,\infty), (-\infty, -x_j+x_i])
\nonumber
\\
&-& G((-x_j,\infty), (-\infty, x_i-x_j])+G((-x_j,\infty), (-\infty, -x_j-x_i])
\nonumber
\\
&+& G((x_j,\infty), (-\infty, x_j+x_i])-G((x_j,\infty), (-\infty, -x_i+x_j])
\nonumber
\\
&=& G((-x_i,x_i], (x_i-x_j,x_j-x_i])
- G((x_i,x_j],(x_j-x_i, x_i+x_j])
\nonumber
\\
&&\qquad
-G((-x_j,-x_i], (-x_i-x_j, x_i-x_j])
\nonumber
\\
&=&\widehat{\Psi} \left(x_i, x_j \right).
\nonumber
\end{eqnarray}
This completes the proof.
\qed

\vskip 10mm

\noindent
{\bf Acknowledgments}
MK would like to thank John Cardy for his hospitality and
useful discussion during his stay in Department of Physics,
Theoretical Physics, University of Oxford, where the present work was done.
HT is partially supported by JSPS Grant-in-Aid for 
Scientific Research Kiban (C) (No. 11640101)
of Japan Society of the Promotion of Science.

\vskip 3mm

\footnotesize 


\vskip 1cm

\noindent
Makoto Katori\\
\\
Department of Physics,\\
Faculty of Science and Engineering,\\
Chuo University, \\
Kasuga, Bunkyo-ku, \\
Tokyo 112-8551, Japan \\
e-mail: katori@phys.chuo-u.ac.jp \\

\vskip 1cm

\noindent
Hideki Tanemura\\
\\
Department of Mathematics and Informatics,\\
Faculty of Science,\\
Chiba University, \\
1-33 Yayoi-cho, Inage-ku,\\
Chiba 263-8522, Japan\\
e-mail: tanemura@math.s.chiba-u.ac.jp\\


\begin{thebibliography}{99}

\bibitem{B72}
B. Belkin,
An invariance principle for conditioned recurrent random walk
attracted to a stable law,
{\it Z.Wahrscheinlichkeitstheorie verw. Geb.} {\bf 21},
(1972), 45-64.

\bibitem{Bi68}
P. Billingsley,
{\it Convergence of Probability Measures},
John Willey \& Sons, 
New York, 1999 (2nd ed.).

\bibitem{Bol76}
E. Bolthausen, 
On a functional central limit theorem for random walks
conditioned to stay positive,
{\it Ann. Probab.} {\bf 4},
(1976), 480-485.

\bibitem{deBr55}
N. G. de Bruijn,
On some multiple integrals involving determinants,
{\it J. Indian Math. Soc.} {\bf 19},
(1955), 133-151.

\bibitem{Do84}
J. L. Doob,
{\it Classical Potential Theory and its Probabilistic Counterpart},
Springer,
1984.

\bibitem{Dur78}
R. Durrett, 
Conditioned limit theorems for some null recurrent Markov processes,
{\it Ann. Probab.} {\bf 6}, 
(1978), 798-828.

\bibitem{Dys62}
F. J. Dyson, 
A Brownian-motion model for the eigenvalues of a random matrix,
{\it J. Math. Phys.} {\bf 3}, 
(1962), 1191-1198.

\bibitem{Fis84}
M. E. Fisher,
Walks, walls, wetting, and melting,
{\it J. Stat. Phys.} {\bf 34}, 
(1984), 667-729.

\bibitem{FH91}
W. Fulton and J. Harris,
{\it Representation Theory},
Springer, 
New York, 1991.

\bibitem{GV85}
I. Gessel and G. Viennot, 
Binomial determinants, paths, and hook length formulae,
{\it Adv. in Math.} {\bf 58}, 
(1985), 300-321.

\bibitem{Gra99}
D. J. Grabiner,
Brownian motion in a Weyl chamber, non-colliding particles,
and random matrices,
{\it Ann. Inst. Henri Poincar\'e} {\bf 35}, 
(1999), 177-204.

\bibitem{Idl74}
D. L. Iglehart,
Functional central limit theorems for random walks
conditioned to stay positive,
{\it Ann. Probab.} {\bf 2}, 
(1974), 608-619.

\bibitem{Imh84}
J. P. Imhof,
Density factorizations for Brownian motion, 
meander and the three-dimensional
Bessel processes, and applications,
{\it J. Appl. Prob.} {\bf 21}, 
(1984), 500-510.

\bibitem{Joh02}
K. Johansson, 
Discrete polynuclear growth and determinantal processes,
math.PR/0206208.

\bibitem{KM59_1}
S. Karlin and L. McGregor,
Coincidence properties of birth and death processes,
{\it Pacific J.} {\bf 9}, 
(1959), 1109-1140.

\bibitem{KM59_2}
S. Karlin and L. McGregor,
Coincidence probabilities,
{\it Pacific J.} {\bf 9}, 
(1959), 1141-1164.

\bibitem{KNT*2}
M. Katori, T. Nagao, and H. Tanemura,
Infinite systems of non-colliding Brownian particles,
to be published in {\it Adv. Stud. in Pure Math.}
``{\it Stochastic Analysis on Large Scale Interacting Systems}",
Mathematical Society of Japan, 2003.

\bibitem{KT02}
M. Katori and H. Tanemura,
Scaling limit of vicious walkers and two-matrix model,
{\it Phys. Rev. E} {\bf 66}, (2002), 011105.

\bibitem{KT03}
M. Katori and H. Tanemura,
Noncolliding Brownian motions and Harish-Chandra formula,
{\it Elect. Comm. in Prob.} {\bf 8}, (2003), 112-121.

\bibitem{KTNK*1}
M. Katori, H. Tanemura, T. Nagao, and
N. Komatsuda,
Vicious walk with a wall, noncolliding meanders,
and chiral and Bogoliubov-deGennes random matrices,
{\it Phys. Rev. E} {\bf 68}, (2003), 021112.

\bibitem{KO01}
W. K\"{o}nig, and N. O'Connell,
Eigenvalues of the Laguerre process 
as non-colliding squared Bessel processes,
{\it Elect. Comm. in Prob.} {\bf 6}, 
(2001), 107-114.

\bibitem{KGV00}
C. Krattenthaler, A. J. Guttmann, and X. G. Viennot, 
Vicious walkers, friendly walkers and Young tableaux: II.
With a wall,
{\it J. Phys. A: Math. Gen.} {\bf 33}, 
(2000), 8835-8866.

\bibitem{Li73}
B. Lindstr\"om, 
On the vector representations of induced matroids,
{\it Bull. London Math. Soc.} {\bf 5}, 
(1973), 85-90.

\bibitem{Mac82}
I. G. Macdonald,
Some conjectures for root systems,
{\it SIAM J. Math. Anal.} {\bf 13}, 
(1982), 988-1007.

\bibitem{Mac95}
I. G. Macdonald,
{\it Symmetric Functions and Hall Polynomials},
Oxford University Press, 
Oxford, 1995 (2nd ed.).

\bibitem{Meh91}
M. L. Mehta, 
{\it Random Matrices}, 
Academic Press, 
London, 1991 (2nd ed.).

\bibitem{MP83}
M. L. Mehta and A. Pandey,
On some Gaussian ensemble of Hermitian matrices,
{\it J. Phys. A: Math. Gen.} {\bf 16}, 
(1983), 2655-2684.

\bibitem{Nagao03}
T. Nagao,
Dynamical correlations for vicious random walk with a wall,
{\it Nucl. Phys.} {\bf B658}[FS], (2003), 373-396.

\bibitem{NF98b}
T. Nagao and P. J. Forrester,
Multilevel dynamical correlation function for Dyson's
Brownian motion model of random matrices,
{\it Phys. Lett.} {\bf A247}, (1998), 42-46.

\bibitem{NKT02}
T. Nagao, M. Katori, and H. Tanemura,
Dynamical correlations among vicious random walkers,
{\it Phys. Lett.} {\bf A307}, (2003), 29-35.

\bibitem{OY02}
N. O'Connell and M. Yor,
A representation for non-colliding random walks,
{\it Elect. Comm. in Prob.} {\bf 7}, 
(2002), 1-12.

\bibitem{Osa96}
H. Osada,
Dirichlet form approach to infinite-dimensional Wiener processes,
{\it Commun. Math. Phys.} {\bf 176}, 
(1996), 117-131.

\bibitem{PM83}
A. Pandey and M. L. Mehta,
Gaussian ensembles of random Hermitian intermediate between
orthogonal and unitary ones,
{\it Commun. Math. Phys.} {\bf 87}, 
(1983), 449-468.

\bibitem{PS02}
M. Pr\"ahofer and H. Spohn, 
Scale invariance of the PNG droplet and the Airy process,
{\it J. Stat. Phys.} {\bf 108}, (2002), 1071-1106.

\bibitem{RY98}
D. Revuz and M. Yor, 
{\it Continuous Martingales and Brownian Motion}, 
Springer, 1998 (3rd ed.).

\bibitem{Shi91}
M. Shimura,
A limit theorem for two dimensional random walk 
conditioned to stay a cone,
{\it Yokohama Math. J.} {\bf 39}, 
(1991), 21-36.

\bibitem{Spiz87}
F. Spitzer,
A Tauberian theorem and its probability interpretation,
{\it Trans. Amer. Math. Soc.} {\bf 94}, 
(1960), 150-169.

\bibitem{Sp87}
H. Spohn,
Interacting Brownian particles:
a study of Dyson's model, in
{\it Hydrodynamic Behavior of Interacting Particle Systems}, 
IMA Volumes in Mathematics and its Applications 9,
ed. G. Papanicolaou, 
Springer, Berlin, 1987.

\bibitem{Stem90}
J. R. Stembridge,
Nonintersecting paths, pfaffians, and the plane partitions,
{\it Adv. in Math.} {\bf 83}, 
(1990), 96-131.

\bibitem{Yor92}
M. Yor,
{\it Some Aspects of Brownian Motion, 
Part I: Some Special Functionals},
Birkh\"auser, Basel, 1992.

\end{thebibliography}
\end{document}